\documentclass[12pt]{article}
\usepackage{graphicx}
\usepackage{color}
\usepackage{amsmath}
\usepackage{amsfonts}
\usepackage{amssymb}
\usepackage{ifthen}
\usepackage{bm,bbm}

\newtheorem{theorem}{Theorem}[section]
\newtheorem{lemma}[theorem]{Lemma}
\newtheorem{corollary}[theorem]{Corollary}

\newtheorem{definition}[theorem]{Definition}

\newcommand{\DEFINED}[1]{{\bf #1}}

\newcommand{\SET}[1]{\mathcal{#1}}

\newcommand{\PRIMEstyle}[1]{{\scriptscriptstyle \mathtt{#1}}}

\newcommand{\itemA}{{\bf a)}}
\newcommand{\itemB}{{\bf b)}}

\newcommand{\Ii}{{\bf i}}

\newcommand{\ZEROvect}{\mathbf{0}}

\newcommand{\ONESvect}{\mathbf{e}}

\newcommand{\STbasis}[1]{e_{#1}}

\newcommand{\REALS}{\mathbbm{R}}

\newcommand{\UNITARY}{\mathcal{U}}

\newcommand{\CONJ}[1]{\overline{#1}}

\newcommand{\ABSOLUTEvalue}[1]{\left| #1 \right|}

\newcommand{\diag}{\mathbf{diag}}

\newcommand{\trace}{\mathbf{tr}}

\renewcommand{\Re}{\mathbf{Re}}
\renewcommand{\Im}{\mathbf{Im}}

\newcommand{\DEFECT}{\mathbf{d}}

\newcommand{\GENdef}{\mathbf{D}}

\newcommand{\FEASIBLEparSPACE}[1]{\mathbbm{V}_{#1}}

\newcommand{\NUMBERofIN}[2]{\sharp_{#2} #1}

\newcommand{\ORTHcomplement}[1]{\left( {#1} \right)^{\bot}}

\renewcommand{\span}{\mathbf{span}}    %original span was already defined

\newcommand{\Mset}[1]{\mathcal{M}_{#1}}

\newcommand{\Msubset}[3]
{\left(\mathcal{M}_{#1}\right)_{\left( {#2} \right), \left( {#3} \right)}}

\newcommand{\Mspace}[1]{\mathbbm{M}_{#1}}

\newcommand{\Msubspace}[3]
{\left(\mathbbm{M}_{#1}\right)_{\left( {#2} \right), \left( {#3}
    \right)}}

\newcommand{\MsubspaceELEMENT}[3]
{{#3}_{\left( {#1} \right), \left( {#2} \right)}}

\newcommand{\MsubspaceDIMENSION}[2]{d^{#1}_{\{#2\}}}

\newcommand{\MsubspaceFIXED}[3]
{\left(\tilde{M}_{#1}\right)_{#2,#3}}

\newcommand{\HADprod}{\circ}

\newcommand{\ELEMENTof}[3]{{\left[ #1 \right]}_{#2,#3}}

\newcommand{\ROWof}[2]{{\left[ #1 \right]}_{#2,:}}

\newcommand{\TRANSPOSE}[1]{{\left( #1 \right)}^T}

\newcommand{\hermTRANSPOSE}[1]{{\left( #1 \right)}^{*}}

\newcommand{\DELTA}[2]
{
  \ifthenelse{\equal{#1}{} \or \equal{#2}{}}
    {\Delta^{#1#2}}
    {\Delta^{#1,#2}}
}

\newcommand{\Mmatrix}[3]{{M}_{#1}^{(#2,#3)}}
\newcommand{\Amatrix}[3]{{A}_{#1}^{(#2,#3)}}
\newcommand{\Smatrix}[3]{{S}_{#1}^{(#2,#3)}}

\newcommand{\ithMATRIX}[2]{{#1}^{(#2)}}

\newcommand{\vectorINDEX}[1]{\left[ #1 \right]}

\newcommand{\PROOFstart}{{\bf Proof}\\}
\newcommand{\PROOFend}{$\blacksquare$}

\begin{document}

\title{Defect of a Kronecker product of unitary matrices}

\author{
  Wojciech Tadej   \\
  \smallskip
  {\small Center for Theoretical Physics, Polish Academy of Sciences, Warsaw, Poland}\\
 \smallskip
 {\small e-mail: wtadej@wp.pl}
%%%% $^1${\small Faculty of Mathematics and Natural Sciences, College of Sciences,}          \\
%%%%     {\small Cardinal Stefan Wyszy{\'n}ski University, Warsaw, Poland}                   \\
%%%%%%%%%  $^1${\small Wydzia\l\ Matematyczno\ -\ Przyrodniczy, Szko\l a Nauk \'Scis\l ych,}         \\
%%%%%%%%%      {\small Uniwersytet Kardyna\l a Stefana Wyszy\'nskiego, Warszawa, Polska (Poland)}    \\
%                                                                                            \\
%   
}
\date{\today}

\maketitle

\begin{abstract}
%The defect $\DEFECT(U)$ of an $N \times N$ unitary matrix $U$ with no
%zero entries is the dimension (called the generalized defect
%$\GENdef(U)$) of the real space of directions, moving into which from
%$U$ we do not disturb the moduli $\ABSOLUTEvalue{U_{i,j}}$ as well as
%the Gram matrix $U^*U$ in the first order, diminished by
%$2N-1$. 
    The generalized defect $\GENdef(U)$ of a unitary $N \times N$
matrix $U$ with no zero entries
is the dimension of the real space of directions, moving into which from
$U$ we do not disturb the moduli $\ABSOLUTEvalue{U_{i,j}}$ as well as
the Gram matrix $U^*U$ in the first order.
    Then the defect $\DEFECT(U)$ is equal to $\GENdef(U)\ -\ (2N-1)$,
that is the generalized defect diminished by the dimension of the manifold
$\{ D_r U D_c\ :\ D_r,D_c\ \mbox{unitary diagonal}\}$. 
    Calculation of $\DEFECT(U)$ involves calculating the dimension
of the space in $\REALS^{N^2}$ spanned by a certain set of vectors
associated with $U$. We split this space into a direct sum, assuming
that $U$ is a Kronecker product of unitary matrices, thus making it
easier to perform calculations numerically. Basing on this, we give a
lower bound on $\GENdef(U)$ (equivalently $\DEFECT(U)$), supposing it
is achieved for most unitaries with a fixed Kronecker product
structure. Also supermultiplicativity of $\GENdef(U)$ with respect to
Kronecker subproducts of $U$ is shown.
\end{abstract}


{\bf MSC-class: 51F25, 47A80, 15B34, 15B51, 58A05}\ \

%@@@@@@@@@@@@@@@@@@@@@@@@@@@@@@@@@@@@@@@@@@@@@@@@@@@@@@@@@@@@@@@@@@@@@@@@@@@@@@@@@@@@@@@@@@@@@@@
%
%                 I  N  T  R  O  D  U  C  T  I  O  N
%
%@@@@@@@@@@@@@@@@@@@@@@@@@@@@@@@@@@@@@@@@@@@@@@@@@@@@@@@@@@@@@@@@@@@@@@@@@@@@@@@@@@@@@@@@@@@@@@@

\section{Introduction}
\label{sec_introduction} The paper is a sequel to \cite{Defect} where the
defect of a unitary matrix was introduced and investigated in a general case.

Let $U \in \UNITARY$, where by $\UNITARY$ we shall denote the set of all
unitaries $N \times N$. Consider the problem: find all directions (complex
matrices) guaranteeing no first order change of $|U_{i,j}|$ for all $i,j$
when moving from $U$ in such a direction, and at the same time no first order
disturbance of unitarity of $U$. The first requirement implies that these
directions are of the form $\Ii R \HADprod U$ where $\HADprod$ is the
Hadamard product and $R$ is a real matrix. The second means that the
directions in question must belong to the space tangent, at $U$, to the $N^2$
dimensional manifold $\UNITARY$. This tangent space is the set $\{EU:\ E\
\mbox{antihermitian}\}$ or, if one likes, $\{UF:\ F\ \mbox{antihermitian}\}$.

Note that the set of directions in question forms a real linear space and it
is overparameterized by $R$ in the case of $U$ having zeros among its
entries. If $U_{i,j}=0$ then the equality $\Ii ( R + \alpha \STbasis{i}
\TRANSPOSE{\STbasis{j}} ) \HADprod U
 \ \ =\ \
 \Ii R \HADprod U$,
with $\STbasis{k}$ being the standard basis column vectors, exposes the fact
that two different $R$'s can be used to parameterize any single direction.
Having this in mind we thus concentrate on the space of $R's$ parameterizing
the space of directions in question:
\begin{equation}
 \label{eq_alowed_directions_parametrizing_space}
 \left\{
   R\ :\ \ \
   \Ii R \HADprod U\ \ \ =\ \ \
   EU\ \ \ \  \mbox{for some antihermitian $E$}
 \right\}
\end{equation}
and ask about its dimension. We define the \DEFINED{defect}\ of $U$ as this
dimension diminished by $2N-1$.

If $U$ has no zero entries, the defect plus $2N-1$ is precisely the dimension
of the space of directions in question. The reduction by $2N-1$ is motivated
by applications of the notion of defect. For most unitaries, especially those
with no zero entries, the $2N$ directions:
\begin{eqnarray}
  \label{eq_phasing_directions}
  \Ii \ \STbasis{i} \TRANSPOSE{\ONESvect} \HADprod U \ \  =\ \
  \Ii \ \diag(\STbasis{i}) \cdot U
  &  &  \mbox{for}\ \ \ i=1..N,\\
%---------------------------------------------------------
  \nonumber
  \Ii \ \ONESvect \TRANSPOSE{\STbasis{j}} \HADprod U \ \  = \ \
  U \cdot \Ii \ \diag(\STbasis{i})
  &  &  \mbox{for}\ \ \ j=1..N,     \\
%---------------------------------------------------------
  \nonumber
  \mbox{where}\ \ \
  \ONESvect\ =\ \TRANSPOSE{[1\ 1\ 1\ \ldots 1]}
  &  &
\end{eqnarray}
which are also directions in question, span the $2N-1$ dimensional real space
tangent to the \DEFINED{phasing manifold}\ $\{ D_r U D_c :\ \ D_r,D_c\
\mbox{unitary diagonal}\} $ generated by left and right multiplication of $U$
by unitary diagonals. The manifold consists of what we regard as unitaries
\DEFINED{equivalent}\ to $U$, with the property of having moduli of entries
identical with those of $U$. 
%In most cases the directions
%(\ref{eq_phasing_directions}) are uninteresting. For more details and
%description of applications of the defect in the search of continuous
%families od unitaries with fixed moduli (among other: families of complex
%Hadamard matricex) we suggest to consult \cite{Defect}.
Usually the directions (\ref{eq_phasing_directions}), leading from $U$
to these equivalent unitaries, are uninteresting.

Other directions in question may lead to unitaries inequivalent to $U$,
but with the same pattern of moduli, provided they exist in a
neighbourhood of $U$. A special case of this property is the existence
of a smooth family of inequivalent complex Hadamard matrices stemming from a given
complex Hadamard matrix $H$. Here by a complex Hadamard matrix we
understand, unlike in combinatorics community, a unitary with all
moduli equal to $1/\sqrt{N}$, although  matrices rescaled to
have unimodular entries are more often considered in this context. As
far as such complex Hadamard matrices are concerned, they have many
applications especially in quantum information theory. We recommend
the reader to consult the online catalogue of complex Hadamard
matrices available at \cite{Catalogue}.

The smallest matrix size for which the known collection of inequivalent complex
Hadamards is presumably incomplete is $6$. Current investigations are
concentrated on the search for a $4$-dimensional smooth family of
inequivalent complex Hadamards. This search is partly motivated by the fact that
the value of the defect calculated for most known $6 \times 6$ complex
Hadamards, including the unitary Fourier matrix $F_6$, is equal to
$4$. This is because the defect of a complex Hadamard $H$ provides an upper bound for the
dimension of a smooth family of inequivalent complex Hadamards
stemming from $H$. We suggest the reader to throw a look into our previous work 
\cite{Defect} to find more details about this application of the
defect. There he can also find that the defect of unitary $U$ (complex
Hadamard $H$) being equal to zero implies that there are no unitaries
(complex Hadamards) inequivalent to $U$ ($H$) in its neighbourhood. At
this point it is worth to recall that analogous (see \cite{Defect}) result was obtained by
Remus Nicoara in his work \cite{Nicoara_span_condition}, under the
name of 'span condition', but in a different context of commuting
squares of algebras, in which complex Hadamard matrices are used to
define commuting squares of matrix algebras.

This above mentioned value $4$, as we now understand, causes one to
suppose that all these $6 \times 6$ complex Hadamards (or complex
Hadamards equivalent to these) with  such a
defect are connected by a single $4$-dimensional family. The old
$2$-dimensional family stemming from the Fourier $F_6$ was reported, 
(among other even earlier sources \cite{Ha96,Di04}) already in our catalogue 
\cite{Catalogue_printed,Catalogue}. Other $1$-dimensional families
were constructed in
\cite{Beauchamp_1dim_selfadjoint,Matolcsi_1dim_symmetric}, 
$2$-dimensional in 
\cite{Szollosi_2dim_hypocycloid,Karlsson_2dim},
$3$ dimensional  in 
\cite{Karlsson_3dim}, 
and the first attempt to build a $4$-dimensional family of $6 \times
6$ inequivalent Hadamards has to be attributed to Ferenc Sz\"{o}ll\H{o}si due
to his succesful work \cite{Szollosi_4dim}, inspired also by some
numerical evidence of the existence of such a family contained in
\cite{Skinner_4dim_numerical}.

From the above collection of results we conclude that it is of value to
study the notion of defect, to search for effective ways to calculate
it or estimate it for unitary matrices of a large size.

In the following we employ another, more expedient characterization of the
defect. To this end let us start with rewriting the definition
(\ref{eq_alowed_directions_parametrizing_space}),
\begin{eqnarray}
 \label{eq_parametrizing_space_equivalent_conditions}
 \Ii R \HADprod U\ \ \ =\ \ \
   EU\ \ \ \  \mbox{for some antihermitian $E$}
 & \Longleftrightarrow &                                            \\
%--------------------------------------------------------------------
 (\Ii R \HADprod U)\hermTRANSPOSE{U}\ \ \ \mbox{is antihermitian}
 & \Longleftrightarrow &  \nonumber                                 \\
%--------------------------------------------------------------------
 U \hermTRANSPOSE{R \HADprod U}\ \ =\ \ (R \HADprod U)
 \hermTRANSPOSE{U}\ \ ,
 &  &  \nonumber
\end{eqnarray}
where $^*$ denotes the Hermitian conjugate.

This further leads to
\begin{eqnarray}
   \forall\  1 \leq i \leq j \leq N
   \ \ \
   \ROWof{U}{i}
   \left[
     \diag(\ROWof{R}{j})
   \right]
   \hermTRANSPOSE{\ROWof{U}{j}}
   & =  &
   \ROWof{U}{i}
   \left[
     \diag(\ROWof{R}{i})
   \right]
   \hermTRANSPOSE{\ROWof{U}{j}}       \nonumber           \\
%------------------------------------------------
   &  \Updownarrow  &      \label{eq_parametrizing_space_equivalent_conditions_entrywise} \\
% -----------------------------------------------
   \forall\  1 \leq i < j \leq N
   \ \ \
   \ROWof{U}{i}
   \left[
     \diag(\ROWof{R}{j})
   \right]
   \hermTRANSPOSE{\ROWof{U}{j}}
   &    =   &
   \ROWof{U}{i}
   \left[
     \diag(\ROWof{R}{i})
   \right]
   \hermTRANSPOSE{\ROWof{U}{j}}\ \ ,    \nonumber
\end{eqnarray}
where by $\ROWof{X}{k}$ the $k$-th row of X is denoted, we
adopt MATLAB notation for rows and columns.

The above can be rewritten as
\begin{equation}
  \label{eq_parametrizing_space_equivalent_conditions_with_trace}
   \forall\ 1 \leq i < j \leq N \ \ \
   \trace{R \TRANSPOSE{\Mmatrix{U}{i}{j}} } \ \ = \ \ 0\ \ ,
\end{equation}
where
\begin{equation}
  \label{eq_ijth_Mmatrix_definition}
  \Mmatrix{U}{i}{j}\ \ =\ \
  \left[
    \begin{array}{cc}
      \cdots\ \ 0\ \ \cdots   &     \\
      %--------------------------------
      \vdots &   \\
      %--------------------------
      \ROWof{U}{i} \HADprod \CONJ{\ROWof{U}{j}}  &  \lefteqn{\ \ \ \scriptstyle{i}}  \\
      %-------------------------------------------
      \vdots &    \\
      %-------------------------------------------
      \cdots\ \ 0\ \ \cdots   &     \\
      %-------------------------------------------
      \vdots  &    \\
      %-------------------------------------------
      -\ \ROWof{U}{i} \HADprod \CONJ{\ROWof{U}{j}}  &  \lefteqn{\ \ \ \scriptstyle{j}} \\
      %-------------------------------------------
      \vdots &    \\
      %-------------------------------------------
      \cdots\ \ 0\ \ \cdots  &    \\
    \end{array}
  \right]
  \ \ = \ \
  \Amatrix{U}{i}{j} \ + \Ii \cdot \Smatrix{U}{i}{j}
\end{equation}
has only its $i$-th and $j$-th rows potentially filled with nonzeros. Here
$\Amatrix{U}{i}{j}$ and $\Smatrix{U}{i}{j}$ denote the real and imaginary
part of $\Mmatrix{U}{i}{j}$.

Let us convert
(\ref{eq_parametrizing_space_equivalent_conditions_with_trace})
into purely real conditions:
\begin{equation}
  \label{eq_parametrizing_space_equivalent_conditions_with_inner_product}
  \forall\ 1 \leq i < j \leq N\ \ \ \ \
  \trace{R \TRANSPOSE{\Amatrix{U}{i}{j}} }\ \ =\ \ 0
  \ \ \ \ \ \mbox{and}\ \ \ \ \
  \trace{R \TRANSPOSE{\Smatrix{U}{i}{j}} }\ \ =\ \ 0\ \ ,
\end{equation}
The form $\trace A\TRANSPOSE{B}$ is a Hilbert-Schmidt type inner product on
the space of real $N \times N$ matrices, hence
(\ref{eq_parametrizing_space_equivalent_conditions_with_inner_product})
amounts to orthogonality conditions in this space. Consequently, and in
accordance with our initial definition the defect of $U$ can be calculated as
\begin{eqnarray}
  \label{eq_defect_calculation_with_orthogonal_complement}
  \dim\left(
    \left\{
      R:\ \Ii R \HADprod U \ \ =\ \ EU\ \ \ \
      \mbox{for some antihermitian $E$}
    \right\}
  \right)
  \ - &(2N-1)& \nonumber
     \\
%---------------------------------------------------
 =\dim\left(
   \ORTHcomplement{
     \span\left(
       \left\{
         \Amatrix{U}{i}{j},\ \Smatrix{U}{i}{j}\
         :\ \ 1 \leq i < j \leq N
       \right\}
     \right)
   }
 \right)
  \ - &(2N-1)&
\end{eqnarray}
where $\ORTHcomplement{}$ is used to indicate the orthogonal complement of
the space  $\Mspace{U}$ spanned by matrices $\Amatrix{U}{i}{j}$,
$\Smatrix{U}{i}{j}$, i.e.,
\begin{equation}
  \label{eq_Mspace_definition}
  \Mspace{U}\ \ \ \ \stackrel{def}{=}\ \ \ \
  \span\left(
    \Mset{U}
  \right)\ ,
\end{equation}
where,
\begin{equation}
  \label{eq_Mset_definition}
  \Mset{U}\ \ \ \ \stackrel{def}{=}\ \ \ \
  \left\{ \Amatrix{U}{i}{j}:\  1 \leq i < j \leq N  \right\}
  \cup
  \left\{ \Smatrix{U}{i}{j}:\ 1 \leq i < j \leq N \right\}\ \ ,
\end{equation}

The properties of orthogonal complements allow now for the following
alternative definition of the defect:

%%%%%%%%%%%%%%%%%%%%%%%%%%%%%%%%%%%%%%%%%%%%%%%%%%%%%%%%%%%%%%%%%%%%%%%%%%
%     D E F I N I T I O N           Defect
%%%%%%%%%%%%%%%%%%%%%%%%%%%%%%%%%%%%%%%%%%%%%%%%%%%%%%%%%%%%%%%%%%%%%%%%%%
\begin{definition}
  \label{def_defect}
  The defect of a unitary $N \times N$ matrix $U$ is the number
  \begin{eqnarray}
    \nonumber
    \DEFECT(U) &  =  &
    \dim\left(
      \ORTHcomplement{
        \Mspace{U}
      }
    \right)\ \ -\ \ (2N-1)         \\
    %------------------------------------
    \label{eq_defect}
    & = &
    (N-1)^2 \ \ -\ \
    \dim\left(
      \Mspace{U}
    \right)\ \ .
  \end{eqnarray}
\end{definition}

It was shown by Karabegov \cite{Karabegov} that the set of those $U$ for
which the defect is greater than $0$ is of the Haar measure zero within the set of all
unitaries $N \times N$. In this article we give a lower bound on the defect
of a Kronecker product of unitary matrices:
\begin{equation}
  \label{eq_kronecker_product_of_unitaries}
  U \ \ =\ \
  \ithMATRIX{U}{1} \otimes \ithMATRIX{U}{2}
  \otimes \ldots \otimes
  \ithMATRIX{U}{r}\ \ ,
\end{equation}
where the unitaries are of size $n_1 \times n_1$, $n_2 \times n_2$, ..., $n_r
\times n_r$,  respectively. 
   The bound computed according to our formulas agrees with the defect
obtained numerically  for Kronecker products of random unitary matrices
drawn according to the Haar measure on the unitary group.

%Our formulas are confirmed numerically to be the
%defect of the Kronecker product of randomly chosen unitary matrices
%of an appropriate size.

However, we do not have here a result concerning measures,  
analogous to that of Karabegov. 
Namely, that for unitary matrices $\ithMATRIX{U}{k}$ lying beyond a
certain set of the Haar measure zero, within the set of all unitary $n_k \times n_k$
matrices, the defect of the resulting Kronecker product 
(\ref{eq_kronecker_product_of_unitaries}) is equal to our bound.

%Namely, that the
%set of Kronecker products of unitaries with the defect greater than our bound is
%infinitely small in the set of all such products with a fixed structure
%(fixed values $n_l$).

Note that the above Kronecker product structure forms a key tool in
studying quantum composite systems. For instance, any local unitary
operation acting on a system composed of $r$ particles has the form 
(\ref{eq_kronecker_product_of_unitaries}). Thus our results on the
defect describe algebraic properties of generic local unitary
operations, which are also called quantum gates.

%@@@@@@@@@@@@@@@@@@@@@@@@@@@@@@@@@@@@@@@@@@@@@@@@@@@@@@@@@@@@@@@@@@@@@@@@@@@@@@@@@@@@@@@@@@@@@@@
%
%                 P R E P A R A T I O N
%
%@@@@@@@@@@@@@@@@@@@@@@@@@@@@@@@@@@@@@@@@@@@@@@@@@@@@@@@@@@@@@@@@@@@@@@@@@@@@@@@@@@@@@@@@@@@@@@@

\section{Preparation}
\label{sec_preparation}

Let us introduce vector indices into Kronecker product
(\ref{eq_kronecker_product_of_unitaries}),
\begin{equation}
  \label{eq_vector_indices}
  \ELEMENTof{U}{\vectorINDEX{i_1..i_r}}{\vectorINDEX{j_1..j_r}}
  \stackrel{def}{=}
  U_{i,j}\ \ ,
\end{equation}
where
\begin{eqnarray}
  \label{eq_vector_index_into_ordinary_index_for_rows}
  i\ \ =\ \ (i_1 - 1) \prod_{k=2}^{r} n_k \ \ +\ \
            (i_2 - 1) \prod_{k=3}^{r} n_k \ \ +\ \
            \ldots \ \
            (i_{r-1} - 1) n_r \ \ +\ \
            i_r
  &  &                                                \\
%---------------------------------------------------------
  j\ \ =\ \ (j_1 - 1) \prod_{k=2}^{r} n_k \ \ +\ \
            (j_2 - 1) \prod_{k=3}^{r} n_k \ \ +\ \
            \ldots \ \
            (j_{r-1} - 1) n_r \ \ +\ \
            j_r
  &  &   \label{eq_vector_index_into_ordinary_index_for_columns}
\end{eqnarray}
and
\begin{equation}
  \label{eq_vectorth_element_calculation}
  \ELEMENTof{U}{\vectorINDEX{i_1..i_r}}{\vectorINDEX{j_1..j_r}}
  \ \ =\ \
  \ELEMENTof{\ithMATRIX{U}{1}}{i_1}{j_1} \ \cdot\
  \ELEMENTof{\ithMATRIX{U}{2}}{i_2}{j_2} \ \cdot\
  \ldots \ \cdot\
  \ELEMENTof{\ithMATRIX{U}{r}}{i_r}{j_r}
\end{equation}
In the following, if not stated otherwise, $i$, $j$ will correspond to the
vector indices $\vectorINDEX{i_1..i_r}$, $\vectorINDEX{j_1..j_r}$.

Let us try for a while to understand vector indices better. The
relation between ordinary and vector indices could be defined in a
different, but completely equivalent way. For two distinct elements in
$\{1..n_1\} \times \ldots \times \{1..n_r\}$ let us introduce a
relation:
\begin{equation}
  \label{eq_vector_indices_inequality}
  \vectorINDEX{i_1..i_r}\ <\ \vectorINDEX{j_1..j_r}
  \ \ \Longleftrightarrow\ \
  i_1=j_1,\ \ldots,\ i_p=j_p,\ i_{p+1} < j_{p+1}
  \ \ \mbox{for some}\ \ p \in \{1,..,r-1\}\ \ .
\end{equation}
Then the relation between ordinary indices and vector indices  we
choose to be the only bijection $\phi$ from\ \
$\{1,\ldots,\prod_{l=1}^r n_l \}$\ \  into\ \
$\{1..n_1\} \times \ldots \times \{1..n_r\}$\ \  (for which we write
$\phi( i ) = \vectorINDEX{i_1..i_r}$), which satisfies:
\ \ $i<j\ \Longrightarrow \phi( i ) < \phi( j )$.

Later in this article we talk about subrows and submatrices, built
from elementary objects (entries of a row for a subrow, rows or
columns of a matrix for a submatrix) indexed by those $i$, whose
$k$-th subindex $i_k$ is fixed at some value:\ \
$i \in \SET{A} = \{i:\ i_k = \alpha\}$. Any $i \in \SET{A}$ is
therefore determined by an element of
$\{1..n_1\} \times \ldots \times \{1..n_{k-1}\} \times
 \{1..n_{k+1}\} \times \ldots \times \{1..n_r\}$, namely the reduced
 vector index $\vectorINDEX{i_1..i_{k-1},i_{k+1}..i_r}$ being a
 subvector of $\vectorINDEX{i_1..i_{k-1},\alpha,i_{k+1}..i_r}$
 corresponding to $i$. We can introduce the relation $<$ for
 the reduced vector indices, and then we note that:
\begin{equation}
  \label{eq_vind_vind_reduced_inequlaity_equivalence}
  \vectorINDEX{i_1..i_{k-1},i_{k+1}..i_r}\ <
  \vectorINDEX{j_1..j_{k-1},j_{k+1}..j_r}
  \ \ \Longleftrightarrow\ \
  \vectorINDEX{i_1..i_{k-1},\alpha,i_{k+1}..i_r}\ <
  \vectorINDEX{j_1..j_{k-1},\alpha,j_{k+1}..j_r}\ \ .
\end{equation}
Ordering the reduced vector indices induces, again, their relation
with reduced ordinary indices
$i',j' \in \{1,\ldots,\prod_{l=1 ,l \neq k}^r n_l \}$ uniquely, that
is it is the same as if we used appropriate versions of formulas
(\ref{eq_vector_index_into_ordinary_index_for_rows}),
(\ref{eq_vector_index_into_ordinary_index_for_columns}).

Let us analyze, as an example, the case of a submatrix $A'$ of matrix
$A$ of size $N \times N$, $N = \prod_{l=1}^r n_l$. Let $A'$ be built
from rows of $A$ indexed by $i \in \SET{A}$, with the order of rows
preserved as is the habit when taking submatrices. Because for $i,j$
indexing rows of $A$ used to build $A'$ and for $i',j'$ corresponding
to reduced vector indices for $i$ and $j$ there holds:
\begin{eqnarray}
  \label{eq_ind_ind_reduced_inequality_equivalence}
  i\ <\ j
  \ \ \Longleftrightarrow\ \
  \vectorINDEX{i_1..i_{k-1},\alpha,i_{k+1}..i_r}
  \ <\ \vectorINDEX{j_1..j_{k-1},\alpha,j_{k+1}..j_r}
  \ \ \Longleftrightarrow\ \
  &  &  \\
%-------------------------------------------------------
  \nonumber
  \vectorINDEX{i_1..i_{k-1},i_{k+1}..i_r}
  \ <\ \vectorINDEX{j_1..j_{k-1},j_{k+1}..j_r}
  \ \ \Longleftrightarrow\ \
  i'\ <\ j'
  &  &   \\
%-------------------------------------------------------
  \nonumber
  \mbox{where $i',j' \in \{1,\ldots, \prod_{l=1,l \neq k}^r n_l\}$}
  &  &
\end{eqnarray}
it is justified to say that the $m$-th row of $A$ forming $A'$ is the
$m'$-th row of $A'$. Here $m'$ is obtained from $m$ in the composition
of relations:
$m\ \longrightarrow\
\vectorINDEX{m_1..m_{k-1},\alpha,m_{k+1}..m_r}\ \longrightarrow\
\vectorINDEX{m_1..m_{k-1},m_{k+1}..m_r}\ \longrightarrow\  m'$.
Equivalently we could say that the
$\vectorINDEX{m_1..m_{k-1},\alpha,m_{k+1}..m_r}$-th row of matrix
$A$ is the
$\vectorINDEX{m_1..m_{k-1},m_{k+1}..m_r}$-th row of its submatrix $A'$.
\bigskip
%-----

Now we return to the main course of this lecture. From now on we
silently assume that all Kronecker factors in
(\ref{eq_kronecker_product_of_unitaries} )  are of size at
least $2$. This we do for convenience of thinking, though we do not
exclude the possibility that some or even all of the below arguments
remain correct without this assumption. We will address this issue
later after announcing the main results at the end of this and at the
beginning of the next section.  

Using the introduced
vector indices, we will split set $\Mset{U}$ defined in
(\ref{eq_Mset_definition}) into disjoint subsets
$\Msubset{U}{k_1,k_2,..,k_p}{v_1,v_2,..,v_p}$,

%%%%%%%%%%%%%%%%%%%%%%%%%%%%%%%%%%%%%%%%%%%%%%%%%%%%%%%%%%
%   D E F I N I T I O N :    M subsets
%%%%%%%%%%%%%%%%%%%%%%%%%%%%%%%%%%%%%%%%%%%%%%%%%%%%%%%%%%

\begin{definition}
  \label{def_Msubsets}
%
  %\begin{itemize}
%    \item
%    \item
          \begin{eqnarray}
            \label{eq_Msubsets_nonempty_parameter}
            \Msubset{U}{k_1,k_2,..,k_p}{v_1,v_2,..,v_p} \ \ =\ \
            \ \ \ \ \ \ \ \ \ \ \ \ \ \ \ \ \ \ \ \ \ \ \ \ \ \ \ \
            \ \ \ \ \ \ \ \ \ \ \ \ \ \ \ \ \ \ \ \ \ \ \ \ \ \ \ \
            \ \ \ \ \ \ \ \ \ \ \ \ \ \ \ \ \ \ \ \ \ \ \ \ \ \ \ \
            &   &              \\
            &   &   \nonumber  \\
          %------------------------------------------------------------
            \left\{
              \Amatrix{U}{i}{j},\ \Smatrix{U}{i}{j}\ :\ \
              \begin{array}{c}
                i_{k_1} = j_{k_1}\  =\  v_1 \in \{1..n_{k_1}\} \\
                i_{k_2} = j_{k_2}\  =\  v_2 \in \{1..n_{k_2}\} \\
                \ldots                                         \\
                i_{k_p} = j_{k_p}\  =\  v_r \in \{1..n_{k_p}\} \\
              \end{array}
              \mbox{and}\ \
              i_k\ \neq\ j_k\ \ \mbox{for}\ \ k \in
              \{1..r\} \setminus \{ k_1,k_2,..,k_p \}
            \right\}\ \ ,
            &   &                                       \nonumber
          \end{eqnarray}
          where we assume that $1 \leq k_1 < \ldots < k_p \leq r$,
          \ \ $1 \leq p \leq r-1$,
          and
          \begin{equation}
            \label{eq_Msubsets_empty_parameter}
            \Msubset{U}{}{}\ \  =\ \
            \left\{
              \Amatrix{U}{i}{j},\ \Smatrix{U}{i}{j}\ :\ \
              i_1 \neq j_1,\
              i_2 \neq j_2,\
              \ldots,\
              i_r \neq j_r
            \right\}\ \ ,
          \end{equation}
          where $()$ means the empty sequence of indices.
%   \end{itemize}
\end{definition}

%In the above we take those $\Amatrix{U}{i}{j}$ and $\Smatrix{U}{i}{j}$
%whose nonzero rows only occure on pairs of positions indicated by
%vector indices with fixed $k_1$-th, $k_2$-th, ..., $k_p$-th entries,
%the other entries differing.

The corresponding subspaces will be denoted according to the following

%%%%%%%%%%%%%%%%%%%%%%%%%%%%%%%%%%%%%%%%%%%%%%%%%%%%%%%%%%%%%%%%%%%%%%%
%  D E F I N I T I O N :   M spaces
%%%%%%%%%%%%%%%%%%%%%%%%%%%%%%%%%%%%%%%%%%%%%%%%%%%%%%%%%%%%%%%%%%%%%%%

\begin{definition}
  \label{def_Msubspaces}
  \begin{eqnarray}
    \label{eq_Msubspaces_empty_parameter}
    \Msubspace{U}{}{} & \stackrel{def}{=} &
    \span\left(  \Msubset{U}{}{} \right)        \\
  %-------------------------------------------------
    \label{eq_Msubspaces_nonempty_parameter}
    \Msubspace{U}{k_1,..,k_p}{v_1,..,v_p} &  \stackrel{def}{=} &
    \span\left(  \Msubset{U}{k_1,..,k_p}{v_1,..,v_p} \right)
  \end{eqnarray}
\end{definition}
Note that in both of the above definitions the notation refers to the
assumed Kronecker product structure of $U$. Later we will meet $U$ deprived
of one or more of its Kronecker factors. Accordingly, the resulting
$U'$ will be assumed to be built from a fewer number of factors, and this
will be reflected in the notation used for subsets of $\Mset{U'}$ and
subspaces of $\Mspace{U'}$.

Any vector in $\Mspace{U}$ (see (\ref{eq_Mspace_definition})) can be written
as a linear combination of components each of which belongs to one of the
above defined subspaces, hence $\Mspace{U}$ is the algebraic sum of
all of them,
\begin{eqnarray}
  \label{eq_Mspace_as_algebraic_sum_of_Msubspaces}
  \lefteqn{\Mspace{U}\ \  =}  &  &  \\
%---------------------------------------------
  \nonumber
  &
  \Msubspace{U}{}{}
  \ \ +\ \
  \sum_{k \in \{1..r\}}
    \sum_{v \in \{1..n_k\}}
      \Msubspace{U}{k}{v}
  \ \ +\ \       &     \\
  %----------------------------------------------------------
  \nonumber
  &
  \ \ +\ \
  \sum_{k_1 < k_2 \in \{1..r\}}
    \sum_{(v_1,v_2) \in \{1..n_{k_1}\} \times \{1..n_{k_2}\}}
      \Msubspace{U}{k_1,k_2}{v_1,v_2}
  \ \ + \ \
  \ldots      &     \\
  %---------------------------------------------------
  \nonumber
  &
  \ldots \ \ + \ \
  \sum_{k_1 < .. < k_{r-1} \in \{1..r\}}
    \sum_{(v_1,..,v_{r-1}) \in \{1..n_{k_1}\} \times .. \times \{1..n_{k_{r-1}}\}}
      \Msubspace{U}{k_1,..,k_{r-1}}{v_1,..,v_{r-1}}
  &
\end{eqnarray}
Our aim is to show that (\ref{eq_Mspace_as_algebraic_sum_of_Msubspaces}) is
in fact a direct sum, which will allow to estimate its dimension. In the
meantime a little preparation is needed.

Using notation similar to that of MATLAB, let, with $y$ at the $l$-th
position,
\begin{eqnarray}
  \label{eq_kronecker_subrow}
  \lefteqn{
    \ELEMENTof{X}{\vectorINDEX{i_1,..,i_r}}
                 {\vectorINDEX{:,..,:,y,:,..,:}}
    \ \ =
  }
  &    &     \\
%-------------------------------------------------------------------------
  \nonumber
  &
  \left[
    \ELEMENTof{X}{\vectorINDEX{i_1,..,i_r}}
                 {\vectorINDEX{1,..,1,y,1,..,1,1}},\ \
    \ELEMENTof{X}{\vectorINDEX{i_1,..,i_r}}
                 {\vectorINDEX{1,..,1,y,1,..,1,2}},\ \
    \ldots,\ \
    \ELEMENTof{X}{\vectorINDEX{i_1,..,i_r}}
                 {\vectorINDEX{n_1,..,n_{l-1},y,n_{l+1},..,n_{r-1},n_r}}
  \right]
  &
\end{eqnarray}
be a subrow of the $\vectorINDEX{i_1,..,i_r}$-th row of matrix $X$ of the
size identical with the size of the Kronecker product
(\ref{eq_kronecker_product_of_unitaries}). The subrow is composed of entries
of $\vectorINDEX{j_1,..,j_r}$-th columns for which $j_l = y$. The horizontal
order of entries is preserved.
%The notation introduced in
%(\ref{eq_kronecker_subrow}) will be used in the lemmas to follow.

%%%%%%%%%%%%%%%%%%%%%%%%%%%%%%%%%%%%%%%%%%%%%%%%%%%%%%%%%%%%%%%%%%%%%%%
%     L E M M A :    Sums of subrows in A, S matrices
%%%%%%%%%%%%%%%%%%%%%%%%%%%%%%%%%%%%%%%%%%%%%%%%%%%%%%%%%%%%%%%%%%%%%%%

\begin{lemma}
  \label{lem_A_S_subrow_properties}
  Let $\vectorINDEX{i_1,..,i_r}$, $\vectorINDEX{j_1,..,j_r}$ be vector
  indices of the Kronecker product (\ref{eq_kronecker_product_of_unitaries})
  corresponding to ordinary indices $i$,$j$.
  \begin{description}
    \item{\itemA}
          If $i_k \neq j_k$ then any $\vectorINDEX{b_1,..,b_r}$-th row
          of $\Amatrix{U}{i}{j}$, $\Smatrix{U}{i}{j}$ satifies:
          \begin{eqnarray}
            \label{eq_A_subrows_zero_sum}
            \sum_{c_k = 1}^{n_k}
              \ELEMENTof{\Amatrix{U}{i}{j}}{\vectorINDEX{b_1,..,b_r}}
                                           {\vectorINDEX{:,..,:,c_k,:,..,:}}
            & = &  \ZEROvect  \\
          %----------------------------------------------------------------
            \label{eq_S_subrows_zero_sum}
            \sum_{c_k = 1}^{n_k}
              \ELEMENTof{\Smatrix{U}{i}{j}}{\vectorINDEX{b_1,..,b_r}}
                                           {\vectorINDEX{:,..,:,c_k,:,..,:}}
            & = &  \ZEROvect
          \end{eqnarray}
          where $c_k$ is at the $k$-th position.

    \item{\itemB} If $i_k = j_k$ then
%the only nonzero rows of $\Amatrix{U}{i}{j}$, $\Smatrix{U}{i}{j}$ have
%this property:
          \begin{eqnarray}
            \label{eq_A_subrows_nonzero_sum}
            \sum_{c_k = 1}^{n_k}
              \ELEMENTof{\Amatrix{U}{i}{j}}{\vectorINDEX{b_1,..,b_r}}
                                           {\vectorINDEX{:,..,:,c_k,:,..,:}}
            & = &
            \ELEMENTof{\Amatrix{U'}{i'}{j'}}{\vectorINDEX{b_1,..,b_{k-1},b_{k+1},..,b_r}}
                                            {:}
            \\
          %------------------------------------------------------------------
            \label{eq_S_subrows_nonzero_sum}
            \sum_{c_k = 1}^{n_k}
              \ELEMENTof{\Smatrix{U}{i}{j}}{\vectorINDEX{b_1,..,b_r}}
                                           {\vectorINDEX{:,..,:,c_k,:,..,:}}
            & = &
            \ELEMENTof{\Smatrix{U'}{i'}{j'}}{\vectorINDEX{b_1,..,b_{k-1},b_{k+1},..,b_r}}
                                            {:}
          \end{eqnarray}
          for any $\vectorINDEX{b_1,..,b_r}$-th row such that
          $b_k = i_k = j_k$, $c_k$ as above.
          On the right hand sides of
          (\ref{eq_A_subrows_nonzero_sum}) and
           (\ref{eq_S_subrows_nonzero_sum}) stand
           $\vectorINDEX{b_1,..,b_{k-1},b_{k+1},..,b_r}$-th rows of
           matrices $\Amatrix{U'}{i'}{j'}$, $\Smatrix{U'}{i'}{j'}$
           constructed from the Kronecker product
           (\ref{eq_kronecker_product_of_unitaries})
           deprived of its $k$-th factor:
          \begin{equation}
            \label{eq_kronecker_product_of_unitaries_reduced}
            U' \ \ =\ \
            \ithMATRIX{U}{1} \otimes \ldots \otimes
            \ithMATRIX{U}{k-1}  \otimes  \ithMATRIX{U}{k+1} \otimes
            \ldots \otimes \ithMATRIX{U}{r}\ \ .
          \end{equation}
          Ordinary indices $i'$, $j'$ correspond to the accordingly
          reduced vector indices \\
          $\vectorINDEX{i_1,..,i_{k-1},i_{k+1},..,i_r}$,
          $\vectorINDEX{j_1,..,j_{k-1},j_{k+1},..,j_r}$.
  \end{description}
\end{lemma}

\PROOFstart     %----------------PROOF  START----------------------------
\itemA\ is of course true for any row of $\Amatrix{U}{i}{j}$,
    $\Smatrix{U}{i}{j}$ indexed by a number different from $i,j$, which is in
    fact a zero row. As the $i$-th and $j$-th rows of $\Amatrix{U}{i}{j}$,
    $\Smatrix{U}{i}{j}$ are of opposite signs, we show \itemA\ only for the
    $i$-th row. We do it by proving that every element of the left hand side
    of (\ref{eq_A_subrows_zero_sum}) and (\ref{eq_S_subrows_zero_sum}) is
    equal to $0$. Indeed for $\Mmatrix{U}{i}{j}$ defined
in  (\ref{eq_ijth_Mmatrix_definition}) we have
\begin{eqnarray}
  \label{eq_M_subrows_zero_sum_entrywise}
  \lefteqn{
    \sum_{c_k = 1}^{n_k}
      \ELEMENTof{\Mmatrix{U}{i}{j}}{\vectorINDEX{i_1,..,i_r}}
                                   {\vectorINDEX{c_1,..,c_r}}
    \ \ =
  }             &     &       \\
%----------------------------------------------------------------
  \nonumber
  &
  \sum_{c_k = 1}^{n_k}
    \left(
      \ELEMENTof{U}{\vectorINDEX{i_1,..,i_r}}
                   {\vectorINDEX{c_1,..,c_r}}
      \cdot
      \CONJ{ \ELEMENTof{U}{\vectorINDEX{j_1,..,j_r}}
                          {\vectorINDEX{c_1,..,c_r}} }
    \right)\ \ =
  &   \\
%----------------------------------------------------------------
  \nonumber
  &
  \sum_{c_k = 1}^{n_k}
    \left(
      \ELEMENTof{\ithMATRIX{U}{1}}{i_1}{c_1}
      \cdot \ldots \cdot
      \ELEMENTof{\ithMATRIX{U}{r}}{i_r}{c_r}
    \right)
    \cdot
    \left(
      \CONJ{\ELEMENTof{\ithMATRIX{U}{1}}{j_1}{c_1}}
      \cdot \ldots \cdot
      \CONJ{\ELEMENTof{\ithMATRIX{U}{r}}{j_r}{c_r}}
    \right)\ \ =
  &     \\
%----------------------------------------------------------------
  \nonumber
  &
  \left(
    \ELEMENTof{\ithMATRIX{U}{k}}{i_k}{1}
    \CONJ{\ELEMENTof{\ithMATRIX{U}{k}}{j_k}{1}}
    \ +\ \ldots\ +\
    \ELEMENTof{\ithMATRIX{U}{k}}{i_k}{n_k}
    \CONJ{\ELEMENTof{\ithMATRIX{U}{k}}{j_k}{n_k}}
  \right)
  \cdot
  \prod_{l=1,\ l \neq k}^{r}
    \left(
      \ELEMENTof{\ithMATRIX{U}{l}}{i_l}{c_l}
      \CONJ{\ELEMENTof{\ithMATRIX{U}{l}}{j_l}{c_l}}
    \right)
  \ \ =\ \
  &              \\
  \nonumber
  &     0\ +\ 0 \Ii\ \ ,     &
\end{eqnarray}
so its real and imaginary parts satisfy (\ref{eq_A_subrows_zero_sum}) and
(\ref{eq_S_subrows_zero_sum}) respectively.

Now part \itemB. If $\vectorINDEX{b_1,..,b_r}$ satisfying $b_k = i_k =
j_k$ is different from
$\vectorINDEX{i_1,..,i_r}$, $\vectorINDEX{j_1,..,j_r}$, then the
ordinary index  $b'$ corresponding to the reduced vector index
$\vectorINDEX{b_1,..,b_{k-1},b_{k+1},..,b_r}$ is different from $i'$,
$j'$. In this situation both sides of
(\ref{eq_A_subrows_nonzero_sum}),
(\ref{eq_S_subrows_nonzero_sum})
are zero rows.

Next we consider the $i$-th row,
$\vectorINDEX{b_1,..,b_r} = \vectorINDEX{i_1,..,i_r}$.
Again, we investigate the sums (\ref{eq_A_subrows_nonzero_sum}) and
(\ref{eq_S_subrows_nonzero_sum}) entrywise. As above, we have:
\begin{eqnarray}
  \label{eq_M_subrows_nonzero_sum_entrywise}
  \lefteqn{
    \sum_{c_k = 1}^{n_k}
      \ELEMENTof{\Mmatrix{U}{i}{j}}{\vectorINDEX{i_1,..,i_r}}
                                   {\vectorINDEX{c_1,..,c_r}}
    \ \ =
  }             &      &     \\
%--------------------------------------------------------------
  \nonumber
  &
  \left(
    \ABSOLUTEvalue{\ELEMENTof{\ithMATRIX{U}{k}}{i_k}{1}}^2
    \ +\ \ldots\ +\
    \ABSOLUTEvalue{\ELEMENTof{\ithMATRIX{U}{k}}{i_k}{n_k}}^2
  \right)
  \cdot
  \prod_{l=1,\ l \neq k}^{r}
    \left(
      \ELEMENTof{\ithMATRIX{U}{l}}{i_l}{c_l}
      \CONJ{\ELEMENTof{\ithMATRIX{U}{l}}{j_l}{c_l}}
    \right)
  \ \ =\ \
  &            \\
%--------------------------------------------------------------
  \nonumber
  &
  \left(
    \prod_{l=1,\ l \neq k}^{r}
      \ELEMENTof{\ithMATRIX{U}{l}}{i_l}{c_l}
  \right)
  \cdot
  \CONJ{
    \left(
      \prod_{l=1,\ l \neq k}^{r}
        \ELEMENTof{\ithMATRIX{U}{l}}{j_l}{c_l}
    \right)
  }\ \ =
  &             \\
%--------------------------------------------------------------
  \nonumber
  &
  \ELEMENTof{U'}{\vectorINDEX{i_1,..,i_{k-1},i_{k+1},..,i_r}}
                {\vectorINDEX{c_1,..,c_{k-1},c_{k+1},..,c_r}}
  \cdot
  \CONJ{
    \ELEMENTof{U'}{\vectorINDEX{j_1,..,j_{k-1},j_{k+1},..,j_r}}
                  {\vectorINDEX{c_1,..,c_{k-1},c_{k+1},..,c_r}}
  }\ \ =
  &               \\
%--------------------------------------------------------------
  \nonumber
  &
  \ELEMENTof{\Mmatrix{U'}{i'}{j'}}{i'}{c'}\ \ ,
  &
\end{eqnarray}
where the ordinary indices $i'$, $j'$, $c'$ correspond to the vectors indices
$\vectorINDEX{i_1,..,i_{k-1},i_{k+1},..,i_r}$,
$\vectorINDEX{j_1,..,j_{k-1},j_{k+1},..,j_r}$,
$\vectorINDEX{c_1,..,c_{k-1},c_{k+1},..,c_r}$ respectively. Now, by taking
the real and imaginary part we get equalities
(\ref{eq_A_subrows_nonzero_sum}) and (\ref{eq_S_subrows_nonzero_sum}).
Calculations for the $j$-th row are completely analogous with the sign
changed.

\PROOFend  %------------PROOF END------------------------------------

The properties stated in the above lemma extend to linear combinations of
matrices $\Amatrix{U}{i}{j},\Smatrix{U}{i}{j}$.

%%%%%%%%%%%%%%%%%%%%%%%%%%%%%%%%%%%%%%%%%%%%%%%%%%%%%%%%%%%%%%%%
%   L E M M A:  a) Subrow of B in M..  is a scaled sum of such subrows
%               b) Subrow of B in M..  add up to a zero row.
%%%%%%%%%%%%%%%%%%%%%%%%%%%%%%%%%%%%%%%%%%%%%%%%%%%%%%%%%%%%%%%%

\begin{lemma}
  \label{lem_B_subrow_properties}
  Let $B \in \Msubspace{U}{k_1,..,k_p}{v_1,..,v_p}$.
  Then:
  \begin{description}
  \item[\itemA] For any $k_s \in \{k_1,..,k_p\}$, any $c_{k_s} \in
      \{1..n_{k_s}\}$, and any  $\vectorINDEX{b_1,..,b_r}$-th row in $B$:
        \begin{equation}
          \label{eq_subrow_as_scaled_sum_of_subrows}
          \ELEMENTof{B}
              {\vectorINDEX{b_1,..,b_r}}
              {\vectorINDEX{:,..,:,c_{k_s},:,..,:}}
          \ \ =\ \
          \ABSOLUTEvalue{
            \ELEMENTof{\ithMATRIX{U}{k_s}}{v_s}{c_{k_s}}
          }^2
          \cdot
          \sum_{d_{k_s} = 1}^{n_{k_s}}
            \ELEMENTof{B}{\vectorINDEX{b_1,..,b_r}}
                         {\vectorINDEX{:,..,:,d_{k_s},:,..,:}}
        \end{equation}
        where $c_{k_s}$, $d_{k_s}$ both at $k_s$-th position designate
        subrows of row $\ELEMENTof{B}{\vectorINDEX{b_1,..,b_r}}{:}$.

        In other words: the $c_{k_s}$-th subrow is a multiple of the sum of
        all subrows designated by all the values of the $k_s$-th
        subindex.

  \item[\itemB]
       For any $k' \notin \{k_1..k_p\}$ and any $b$-th row in $B$:
       \begin{equation}
         \label{eq_subrows_add_up_to_zero_row}
         \sum_{c_{k'}=1}^{n_{k'}}
           \ELEMENTof{B}{b}{\vectorINDEX{:,..,:,c_{k'},:,..,:}}
         =
         \ZEROvect\ \ \ \ \ \mbox{where $c_{k'}$ is at $k'$-th position.}
       \end{equation}
  \end{description}
\end{lemma}

\PROOFstart    %-------PROOF  START--------------------------------

The part \itemA\ of the lemma is of course true for all rows indexed by
$\vectorINDEX{b_1,..,b_r}$ such that $b_{k_1} \neq v_1$ or $b_{k_2} \neq v_2$
or ... or $b_{k_p} \neq v_p$, since these are zero rows coming from zero rows
in $\Amatrix{U}{i}{j}$, $\Smatrix{U}{i}{j}$ spanning
$\Msubspace{U}{k_1,..,k_p}{v_1,..,v_p}$.

Let us thus consider the $\vectorINDEX{b_1,..,b_r}$-th row of B with $b_{k_1}
= v_1$, $b_{k_2} = v_2$, ..., $b_{k_p} = v_p$. Since $B \in
\Msubspace{U}{k_1,..,k_p}{v_1,..,v_p}$ for the $\vectorINDEX{b_1,..,b_r}$-th
row we have,
\begin{eqnarray}
  \label{eq_B_row_as_combination_of_A_S_rows}
  \lefteqn{
    \ELEMENTof{B}{\vectorINDEX{b_1,..,b_r}}{:}\ \ =
  }  &   &    \\
%-------------------------------------------------------
  \nonumber
  &
  \sum_{\Amatrix{U}{i}{j} \in \Msubset{U}{k_1,..,k_p}{v_1,..,v_p}}
    \alpha^{(i,j)}
    \ELEMENTof{\Amatrix{U}{i}{j}}{\vectorINDEX{b_1,..,b_r}}{:}
  \ \ \ +   &     \\
%------------------------------------------------------------------
  \nonumber
  &
  \sum_{\Smatrix{U}{i}{j} \in \Msubset{U}{k_1,..,k_p}{v_1,..,v_p}}
    \sigma^{(i,j)}
    \ELEMENTof{\Smatrix{U}{i}{j}}{\vectorINDEX{b_1,..,b_r}}{:}
  \ \ ,     &
\end{eqnarray}
and analogously for the $c_{k_s}$-th subrow designated by the value of
the $k_s$-th index, where we ommit the ranges for $\alpha^{(i,j)}$ and
$\sigma^{(i,j)}$:
\begin{eqnarray}
  \label{eq_B_subrow_as_combination_of_A_S_subrows}
  \lefteqn{
    \ELEMENTof{B}{\vectorINDEX{b_1,..,b_r}}
                 {\vectorINDEX{:,..,:,c_{k_s},:,..,:}}\ \ =
  }   &    &    \\
%--------------------------------------------------------------
  \nonumber
  &
  \sum_{....}%\Amatrix{U}{i}{j} \in \Msubset{U}{k_1,..,k_p}{v_1,..,v_p}}
    \alpha^{(i,j)}
    \ELEMENTof{\Amatrix{U}{i}{j}}{\vectorINDEX{b_1,..,b_r}}
                                 {\vectorINDEX{:,..,:,c_{k_s},:,..,:}}
  \ \ \ +\ \ \
  \sum_{....}    %\Smatrix{U}{i}{j} \in \Msubset{U}{k_1,..,k_p}{v_1,..,v_p}}
    \sigma^{(i,j)}
    \ELEMENTof{\Smatrix{U}{i}{j}}{\vectorINDEX{b_1,..,b_r}}
                                 {\vectorINDEX{:,..,:,c_{k_s},:,..,:}}
  \ \ .  &
\end{eqnarray}

Hence, on the one hand:
\begin{eqnarray}
  \label{eq_sum_of_subrows_in_linear_combination}
  \lefteqn{
    \sum_{d_{k_s}=1}^{n_{k_s}}
      \ELEMENTof{B}{\vectorINDEX{b_1,..,b_r}}
                   {\vectorINDEX{:,..,:,d_{k_s},:,..,:}}
    \ \ =
  }          &      &     \\
%-------------------------------------------------------
  \nonumber
  &
  \sum_{....}
    \alpha^{(i,j)}
    \sum_{d_{k_s}=1}^{n_{k_s}}
      \ELEMENTof{\Amatrix{U}{i}{j}}{\vectorINDEX{b_1,..,b_r}}
                                   {\vectorINDEX{:,..,:,d_{k_s},:,..,:}}
  \ \ +\ \
  \sum_{....}
    \sigma^{(i,j)}
    \sum_{d_{k_s}=1}^{n_{k_s}}
      \ELEMENTof{\Smatrix{U}{i}{j}}{\vectorINDEX{b_1,..,b_r}}
                                   {\vectorINDEX{:,..,:,d_{k_s},:,..,:}}
  \ \ =   &    \\
%-----------------------------------------------------------
  \nonumber
  &
  \sum_{....}
    \alpha^{(i,j)}
    \ELEMENTof{\Amatrix{U'}{i'}{j'}}
              {\vectorINDEX{b_1,..,b_{k_s-1},b_{k_s+1},..,b_r}}
              {:}
  \ \ +\ \
  \sum_{....}
    \sigma^{(i,j)}
    \ELEMENTof{\Smatrix{U'}{i'}{j'}}
              {\vectorINDEX{b_1,..,b_{k_s-1},b_{k_s+1},..,b_r}}
              {:}
  &
\end{eqnarray}
where we used Lemma \ref{lem_A_S_subrow_properties} \itemB, which can be
employed since the indices $i$,$j$ in the sums satisfy $i_{k_1} = j_{k_1} =
v_1$, ..., $i_{k_p} = j_{k_p} = v_p$, and in particular $i_{k_s} =
j_{k_s}$, as well as $b_{k_s} = v_s = i_{k_s} = j_{k_s}$.
As previously, the ordinary indices $i'$, $j'$ correspond to vector indices
$\vectorINDEX{i_1,..,i_{k_s-1},i_{k_s+1},...,i_r}$,
$\vectorINDEX{j_1,..,j_{k_s-1},j_{k_s+1},...,j_r}$ and designate nonzero rows
in matrices $\Amatrix{U'}{i'}{j'}$, $\Smatrix{U'}{i'}{j'}$ constructed from
the reduced Kronecker product $\ithMATRIX{U}{1} \otimes .. \otimes
\ithMATRIX{U}{k_s-1} \otimes \ithMATRIX{U}{k_s+1} \otimes ..
\otimes \ithMATRIX{U}{r}$.

On the other hand  the expression
(\ref{eq_B_subrow_as_combination_of_A_S_subrows}) translates into
\begin{eqnarray}
  \label{eq_B_subrow_as_combination_of_A_S_subrows_continued}
  \lefteqn{
    \ELEMENTof{B}{\vectorINDEX{b_1,..,b_r}}
                 {\vectorINDEX{:,..,:,c_{k_s},:,..,:}}\ \ =
  }   &    &    \\
%-----------------------------------------------------------
  \nonumber
  &
  \sum_{....}
    \alpha^{(i,j)}
    \ABSOLUTEvalue{
      \ELEMENTof{\ithMATRIX{U}{k_s}}{v_s}{c_{k_s}}
    }^2
    \ELEMENTof{\Amatrix{U'}{i'}{j'}}
              {\vectorINDEX{b_1,..,b_{k_s-1},b_{k_s+1},..,b_r}}
              {:}
  \ \ +\ \  &  \\
%-----------------------------------------------------------
  \nonumber
  &
  \sum_{....}
    \sigma^{(i,j)}
    \ABSOLUTEvalue{
      \ELEMENTof{\ithMATRIX{U}{k_s}}{v_s}{c_{k_s}}
    }^2
    \ELEMENTof{\Smatrix{U'}{i'}{j'}}
              {\vectorINDEX{b_1,..,b_{k_s-1},b_{k_s+1},..,b_r}}
              {:}\ \ ,
  &
\end{eqnarray}
with $U'$, $i'$, $j'$ described above. Combining
(\ref{eq_sum_of_subrows_in_linear_combination}) and
(\ref{eq_B_subrow_as_combination_of_A_S_subrows_continued}) we obtain \itemA.

%.........................................PART B
To prove the part {\itemB}  we use again
(\ref{eq_B_row_as_combination_of_A_S_rows}) and rewrite it for the subrows
corresponding to $k' \notin \{k_1..k_p\}$:
\begin{eqnarray}
  \label{eq_B_subrow_as_combination_of_A_S_subrows_with_k_prim}
  \lefteqn{
    \ELEMENTof{B}{\vectorINDEX{b_1,..,b_r}}
                 {\vectorINDEX{:,..,:,c_{k'},:,..,:}}\ \ =
  }   &    &  \nonumber  \\
%--------------------------------------------------------------
  &
  \sum_{....}%\Amatrix{U}{i}{j} \in \Msubset{U}{k_1,..,k_p}{v_1,..,v_p}}
    \alpha^{(i,j)}
    \ELEMENTof{\Amatrix{U}{i}{j}}{\vectorINDEX{b_1,..,b_r}}
                                 {\vectorINDEX{:,..,:,c_{k'},:,..,:}}
  \ \ \ +\ \ \
  \sum_{....}    %\Smatrix{U}{i}{j} \in \Msubset{U}{k_1,..,k_p}{v_1,..,v_p}}
    \sigma^{(i,j)}
    \ELEMENTof{\Smatrix{U}{i}{j}}{\vectorINDEX{b_1,..,b_r}}
                                 {\vectorINDEX{:,..,:,c_{k'},:,..,:}}
\end{eqnarray}
Here $i,j$ are such that $i_{k'} \neq j_{k'}$. We may thus use Lemma
\ref{lem_A_S_subrow_properties} \itemA\ to find that summation of
subrows (\ref{eq_B_subrow_as_combination_of_A_S_subrows_with_k_prim})
produces a zero vector:
\begin{eqnarray}
  \label{eq_sum_of_subrows_in_linear_combination_with_k_prim}
  \lefteqn{
    \sum_{c_{k'}=1}^{n_{k'}}
      \ELEMENTof{B}{\vectorINDEX{b_1,..,b_r}}
                   {\vectorINDEX{:,..,:,c_{k'},:,..,:}}
    \ \ =
  }          &      &     \\
%-------------------------------------------------------
  \nonumber
  &
  \sum_{....}
    \alpha^{(i,j)}
    \sum_{c_{k'}=1}^{n_{k'}}
      \ELEMENTof{\Amatrix{U}{i}{j}}{\vectorINDEX{b_1,..,b_r}}
                                   {\vectorINDEX{:,..,:,c_{k'},:,..,:}}
  \ \ +\ \
  \sum_{....}
    \sigma^{(i,j)}
    \sum_{c_{k'}=1}^{n_{k'}}
      \ELEMENTof{\Smatrix{U}{i}{j}}{\vectorINDEX{b_1,..,b_r}}
                                   {\vectorINDEX{:,..,:,c_{k'},:,..,:}}
  \ \ =   &    \\
%-----------------------------------------------------------
  \nonumber
  &
  \sum_{....}
    \alpha^{(i,j)}
    \cdot \ZEROvect
  \ \ +\ \
  \sum_{....}
    \sigma^{(i,j)}
    \cdot \ZEROvect \ \ \ =\ \ \  \ZEROvect\ \ .
  &
\end{eqnarray}
\PROOFend   %------------------PROOF END-----------------------

Related is the following result which we will also need to prove the
theorem on a direct sum that comes next.

%%%%%%%%%%%%%%%%%%%%%%%%%%%%%%%%%%%%%%%%%%%%%%%%%%%%%%%%%%%%%%%%%%%%
%   L E M M A :      Isomorphism between M-subspaces
%%%%%%%%%%%%%%%%%%%%%%%%%%%%%%%%%%%%%%%%%%%%%%%%%%%%%%%%%%%%%%%%%%%%

\begin{lemma}
  \label{lem_isomorphic_Msubspaces}
  Let $\MsubspaceFIXED{U}{k}{v}$ be the space spanned by matrices
  $\Amatrix{U}{i}{j}$, $\Smatrix{U}{i}{j}$ such that $i_k = j_k = v$,
  where $v \in \{1..n_k\}$. There exists an isomorphism mapping
  $\MsubspaceFIXED{U}{k}{v}$ onto $\Mspace{U'}$,
  where
  $U' = \ithMATRIX{U}{1} \otimes \ldots \otimes
        \ithMATRIX{U}{k-1} \otimes \ithMATRIX{U}{k+1} \otimes \ldots \otimes
        \ithMATRIX{U}{r}$. It can be chosen in such a way that any subspace
  $\Msubspace{U}{k_1,..,k_p}{v_1,..,v_p}$ of
  $\MsubspaceFIXED{U}{k}{v}$ with $k_q = k$ and $v_q = v$ for some
$q \in \{1..p\}$, is mapped onto the subspace
  $\Msubspace{U'}{k'_1,..,k'_{p-1}}
                 {v_1,..,v_{q-1},v_{q+1},..v_p}$
  of $\Mspace{U'}$, where $k'_1$, ..., $k'_{p-1}$ are positions of
  $k_1$, ..., $k_{q-1}$, $k_{q+1}$, ..., $k_{p}$ in the sequence
  $(1, 2, ..., k-1, k+1, ..., r)$. (Note that $k'_g$ indicate positions of factors
  $\ithMATRIX{U}{k_1}$, ..., $\ithMATRIX{U}{k_{q-1}}$,
  $\ithMATRIX{U}{k_{q+1}}$, ..., $\ithMATRIX{U}{k_p}$
  in the product $U'$.) In particular the chosen isomorphism maps
  $\Msubspace{U}{k}{v}$ onto $\Msubspace{U'}{}{}$.
\end{lemma}

\PROOFstart %----------------------------------------------------------
We will show that matrices $\Amatrix{U}{i}{j}$, $\Smatrix{U}{i}{j}$ spanning
$\MsubspaceFIXED{U}{k}{v}$ remain in a simple relation with matrices
$\Amatrix{U'}{i'}{j'}$, $\Smatrix{U'}{i'}{j'}$. Here $i'$, $j'$ correspond to
vector indices $\vectorINDEX{i_1,..,i_{k-1},i_{k+1},..,i_r}$,
$\vectorINDEX{j_1,..,j_{k-1},j_{k+1},..,j_r}$ of the reduced Kronecker
product $U'$ (\ref{eq_kronecker_product_of_unitaries_reduced}), i.e., are
given by the appropriate versions  of
(\ref{eq_vector_index_into_ordinary_index_for_rows}) and
(\ref{eq_vector_index_into_ordinary_index_for_columns}).

We will be more precise about the above mentioned relation later, but at first let
us state it without too many details. Any $\Amatrix{U}{i}{j}$ spanning
$\MsubspaceFIXED{U}{k}{v}$ contains, as its only nonzero and disjoint
submatrices, $n_k$ copies of the corresponding $\Amatrix{U'}{i'}{j'}$ , where
each copy is multiplied by one of $n_k$ multipliers. Column positions of the
copies do not depend on $i,j$. The multipliers, of which at least one is
nonzero, depend only on the column positions of the copies they act on. The
same applies to $\Smatrix{U}{i}{j}$ and $\Smatrix{U'}{i'}{j'}$, with the
pattern of submatrices and values of the multipliers identical with that for
$\Amatrix{U}{i}{j}$ and $\Amatrix{U'}{i'}{j'}$.

For a more detailed explanation let us consider the $i$-th row of
$\Amatrix{U}{i}{j}$. The \\
$\vectorINDEX{d_1,..,d_{k-1},d_{k+1},..,d_r}$-th
element of the $d_k$-th subrow of it, $d_k \in \{1..n_k\}$, reads:
\begin{eqnarray}
  \label{eq_subrow_of_A_row}
  \lefteqn{
    \ELEMENTof{\Amatrix{U}{i}{j}}
              {i}
              {\vectorINDEX{d_1,..,d_k,..,d_r}}
     =
    \Re\left(
      \ELEMENTof{\Mmatrix{U}{i}{j}}
                {i}
                {\vectorINDEX{d_1,..,d_k,..,d_r}}
    \right)
    =
  }
  &   &     \\
%---------------------------------------------------------------------
  \nonumber
  &
  \Re\left(
    \ELEMENTof{U}
              {i}
              {\vectorINDEX{d_1,..,d_k,..,d_r}}
    \cdot
    \ELEMENTof{\CONJ{U}}
              {j}
              {\vectorINDEX{d_1,..,d_k,..,d_r}}
  \right)\ \ ,
  &
\end{eqnarray}
which is equal to
\begin{equation}
  \label{eq_subrow_of_A_row_continuation_1}
  \Re\left(
    \ELEMENTof{\ithMATRIX{U}{1}}{i_1}{d_1}
    \ELEMENTof{\CONJ{\ithMATRIX{U}{1}}}{j_1}{d_1}
    \cdot \ldots \cdot
    \ELEMENTof{\ithMATRIX{U}{k}}{i_k}{d_k}
    \ELEMENTof{\CONJ{\ithMATRIX{U}{k}}}{j_k}{d_k}
    \cdot \ldots \cdot
    \ELEMENTof{\ithMATRIX{U}{r}}{i_r}{d_r}
    \ELEMENTof{\CONJ{\ithMATRIX{U}{r}}}{j_r}{d_r}
  \right)\ \ ,
\end{equation}
where
$\ELEMENTof{\ithMATRIX{U}{k}}{i_k}{d_k}
 \ELEMENTof{\CONJ{\ithMATRIX{U}{k}}}{j_k}{d_k}$
can also be written as
$\ELEMENTof{\ithMATRIX{U}{k}}{v}{d_k}
    \ELEMENTof{\CONJ{\ithMATRIX{U}{k}}}{v}{d_k}
\ =\
\ABSOLUTEvalue{\ELEMENTof{\ithMATRIX{U}{k}}{v}{d_k}}^2$.
As a result we get
\begin{eqnarray}
  \label{eq_subrow_of_A_row_continuation_2}
  \lefteqn{}
  &
  \ABSOLUTEvalue{\ELEMENTof{\ithMATRIX{U}{k}}{v}{d_k}}^2
  \cdot
  \Re\left(
    \prod_{s \in \{1..r\} \setminus {k}}
      \ELEMENTof{\ithMATRIX{U}{s}}{i_s}{d_s}
      \ELEMENTof{\CONJ{\ithMATRIX{U}{s}}}{j_s}{d_s}
  \right)\ \  =
  &     \\
%--------------------------------------------------------
  \nonumber
  &
  \ABSOLUTEvalue{\ELEMENTof{\ithMATRIX{U}{k}}{v}{d_k}}^2
  \cdot
  \Re\left(
    \ELEMENTof{U'}
              {\vectorINDEX{i_1,..,i_{k-1},i_{k+1},..,i_r}}
              {\vectorINDEX{d_1,..,d_{k-1},d_{k+1},..,d_r}}
    \cdot
    \ELEMENTof{\CONJ{U'}}
              {\vectorINDEX{j_1,..,j_{k-1},j_{k+1},..,j_r}}
              {\vectorINDEX{d_1,..,d_{k-1},d_{k+1},..,d_r}}
  \right)\ \ =
  &    \\
%--------------------------------------------------------
  \nonumber
  &
  \ABSOLUTEvalue{\ELEMENTof{\ithMATRIX{U}{k}}{v}{d_k}}^2
  \cdot
  \Re\left(
    \ELEMENTof{\Mmatrix{U'}{i'}{j'}}
              {i'}
              {\vectorINDEX{d_1,..,d_{k-1},d_{k+1},..,d_r}}
  \right)
  \ \ =\ \
  \ABSOLUTEvalue{\ELEMENTof{\ithMATRIX{U}{k}}{v}{d_k}}^2
  \cdot
  \ELEMENTof{\Amatrix{U'}{i'}{j'}}
            {i'}
            {\vectorINDEX{d_1,..,d_{k-1},d_{k+1},..,d_r}}\ \ .
  &
\end{eqnarray}
The equality of the far left hand side of (\ref{eq_subrow_of_A_row})
and the far right hand side of (\ref{eq_subrow_of_A_row_continuation_2})
is
an entry by entry equality between two rows. Now, if we change the
sign of $\Re$ in the expressions in (\ref{eq_subrow_of_A_row}),
(\ref{eq_subrow_of_A_row_continuation_1}) and
(\ref{eq_subrow_of_A_row_continuation_2}) we will obtain a relation
between
  the $d_k$-th subrow of the $j$-th row in $\Amatrix{U}{i}{j}$ and
  the $j'$-th row of $\Amatrix{U'}{i'}{j'}$.
If we used $\Im$ instead of $\Re$, we would get a
relation between
  the $d_k$-th subrow of the $i$-th or $j$-th row of $\Smatrix{U}{i}{j}$ and, correspondingly,
  the $i'$-th or $j'$-th row of $\Smatrix{U'}{i'}{j'}$.
This is summarized below:
\begin{eqnarray}
  \label{eq_AS_U_subrow_AS_U_prim_row_relations}
  \ELEMENTof{\Amatrix{U}{i}{j}}
            {i}
            {\vectorINDEX{:,..,:,d_k,:,..,:}}
  \ \ \ =\ \ \
  \mu_k
  \ \ \cdot\ \
  \ELEMENTof{\Amatrix{U'}{i'}{j'}}
            {i'}
            {:}\ \ ,
  &  &  \\
%----------------------------------------------------------------
  \nonumber
  \ELEMENTof{\Amatrix{U}{i}{j}}
            {j}
            {\vectorINDEX{:,..,:,d_k,:,..,:}}
  \ \ \ =\ \ \
  \mu_k
  \ \ \cdot\ \
  \ELEMENTof{\Amatrix{U'}{i'}{j'}}
            {j'}
            {:}\ \ ,
  &  &  \\
%---------------------------------------------------------------
  \nonumber
  \ELEMENTof{\Smatrix{U}{i}{j}}
            {i}
            {\vectorINDEX{:,..,:,d_k,:,..,:}}
  \ \ \ =\ \ \
  \mu_k
  \ \ \cdot\ \
  \ELEMENTof{\Smatrix{U'}{i'}{j'}}
            {i'}
            {:}\ \ ,
  &  &  \\
%---------------------------------------------------------------
  \nonumber
  \ELEMENTof{\Smatrix{U}{i}{j}}
            {j}
            {\vectorINDEX{:,..,:,d_k,:,..,:}}
  \ \ \ =\ \ \
  \mu_k
  \ \ \cdot\ \
  \ELEMENTof{\Smatrix{U'}{i'}{j'}}
            {j'}
            {:}\ \ ,
  &  &
\end{eqnarray}
where $\mu_k = \ABSOLUTEvalue{\ELEMENTof{\ithMATRIX{U}{k}}{v}{d_k}}^2$.
One of the $n_k$ numbers $\mu_k$ is nonzero since $\ithMATRIX{U}{k}$ is unitary.

Now let us use
$\ELEMENTof{X}{\vectorINDEX{:,..,:,b_k,:,..,:}}
              {\vectorINDEX{:,..,:,c_k,:,..,:}}$
to denote a submatrix of $X$ built of those $x$-th rows and $y$-th columns of
X  for which $x_k = b_k$ and $y_k = c_k$. As a consequence of relations
(\ref{eq_AS_U_subrow_AS_U_prim_row_relations}) we see that
$\Amatrix{U}{i}{j}$, $\Smatrix{U}{i}{j}$ spanning $\MsubspaceFIXED{U}{k}{v}$
are all zero matrices except for their submatrices:
\begin{eqnarray}
  \label{eq_A_nonzero_submatrix}
  \ELEMENTof{\Amatrix{U}{i}{j}}{\vectorINDEX{:,..,:,v,:,..,:}}
                               {\vectorINDEX{:,..,:,c_k,:,..,:}}
  \ \ =\ \
  \mu_k \cdot \Amatrix{U'}{i'}{j'}\ \ \ \ \ c_k \in \{1,..,n_k\}
  \ \ ,
  &  &  \\
%-------------------------------------------------------------
  \label{eq_S_nonzero_submatrix}
  \ELEMENTof{\Smatrix{U}{i}{j}}{\vectorINDEX{:,..,:,v,:,..,:}}
                               {\vectorINDEX{:,..,:,c_k,:,..,:}}
  \ \ =\ \
  \mu_k \cdot \Smatrix{U'}{i'}{j'}\ \ \ \ \ c_k \in \{1,..,n_k\}
  \ \ ,
  &  &      \\
%-------------------------------------------------------------
  \nonumber
  \mbox{where}\ \ \
  \mu_k\ =\ \ABSOLUTEvalue{\ELEMENTof{\ithMATRIX{U}{k}}{v}{c_k}}^2
  &  &
\end{eqnarray}
(where $v$ is at $k$-th position), with one of $\mu_k$ being nonzero.
Analogous equalities we could write for linear combinations on both sides.
( (\ref{eq_A_nonzero_submatrix}) and (\ref{eq_S_nonzero_submatrix})
can be checked entrywise, in which one uses that the
$\vectorINDEX{b_1,..,b_{k-1},b_{k+1},..,b_r},
 \vectorINDEX{c_1,..,c_{k-1},c_{k+1},..,c_r}$-th
element of
$\ELEMENTof{X}{\vectorINDEX{:,..,:,b_k,:,..,:}}
              {\vectorINDEX{:,..,:,c_k,:,..,:}}$
is equal to
$\ELEMENTof{X}{\vectorINDEX{b_1,..,b_r}}
              {\vectorINDEX{c_1,..,c_r}}$.
Here of course we identify, for example,
$\vectorINDEX{b_1,..,b_{k-1},b_{k+1},..,b_r}$ with ordinary index $b'$.
See also the comment that follows the introduction of vector indices.)

Since all $\mu_k$ add up to 1 , i.e. the norm of the $v$-th row of
$\ithMATRIX{U}{k}$, we can choose an isomorphism $\Psi$ between
$\MsubspaceFIXED{U}{k}{v}$ and $\Mspace{U'}$ to be defined as (again assuming
that $v$ stands at $k$-th position):
\begin{equation}
  \label{eq_Psi_isomorphism}
  \Psi\left(  X \in \MsubspaceFIXED{U}{k}{v}  \right)
  \ \ =\ \
  \sum_{c_k = 1}^{n_k}
    \ELEMENTof{X}{\vectorINDEX{:,..,:,v,:,..,:}}
                 {\vectorINDEX{:,..,:,c_k,:,..,:}}\ \ .
\end{equation}
so that
\begin{eqnarray}
  \label{eq_Psi_mapping_AS_matrices}
  \Psi\left( \Amatrix{U}{i}{j} \right) \ \  =\ \  \Amatrix{U'}{i'}{j'}
  &  &   \\
%-----------------------------------------------------------
  \nonumber
  \Psi\left( \Smatrix{U}{i}{j} \right) \ \  =\ \  \Smatrix{U'}{i'}{j'}
\end{eqnarray}

Consider now $\Msubspace{U}{k_1,..,k_p}{v_1,..,v_p} \subset
 \MsubspaceFIXED{U}{k}{v}$ with $k_q = k$ and $v_q = v$ for some $q
 \in \{1,..,p\}$.
Matrices $\Amatrix{U}{i}{j}$,
$\Smatrix{U}{i}{j}$ spanning this space are indexed by $i,j$ for which
$i_{k_1} = j_{k_1} = v_1$, ..., $i_{k_q} = j_{k_q} = v$, ..., $i_{k_p} =
j_{k_p} = v_p$, and finally $i_l \neq j_l$ for $l \notin \{k_1,..,k_p\}$. Note
that according to (\ref{eq_Psi_mapping_AS_matrices}) the isomorphism $\Psi$
maps these matrices into all matrices spanning
$\Msubspace{U'}{k'_1,..,k'_{p-1}}
               {v_1,..,v_{q-1},v_{q+1},..,v_p}$,
where $k'_1$, ..., $k'_{p-1}$ are defined in the lemma.
This is because
if we take any pair $i,j$ designating
   $\Amatrix{U}{i}{j}$, $\Smatrix{U}{i}{j}$
spanning
   $\Msubspace{U}{k_1,..,k_p}{v_1,..,v_p}$,
in the reduced  vector index $\vectorINDEX{i_1,..,i_{k-1},i_{k+1},..,i_r}$ at
positions
    $k'_1$, ..., $k'_{p-1}$
we have subindices
    $i_{k_1}$, ..., $i_{k_{q-1}}$, $i_{k_{q+1}}$, ... $i_{k_p}$
that are equal in pairs with subindices
    $j_{k_1}$, ..., $j_{k_{q-1}}$, $j_{k_{q+1}}$, ... $j_{k_p}$
at the same positions in the vector index
$\vectorINDEX{j_1,..,j_{k-1},j_{k+1},..,j_r}$. Their values are $v_1$, ...,
$v_{q-1}$, $v_{q+1}$, ..., $v_{p}$, respectively. The only other pair of
equal subindices, namely $i_{k_q} = j_{k_q} = v$, is absent from the reduced
vector indices corresponding to $i',j'$. At other positions we get any
possible pair of different values as we take any pair $i,j$ allowed here,
that is indexing $\Amatrix{U}{i}{j}$, $\Smatrix{U}{i}{j}$ spanning
$\Msubspace{U}{k_1,..,k_p}{v_1,..,v_p}$.
%REPETITION
%Let us repeat: by taking, for example, all $\Amatrix{U}{i}{j}$ spanning
%$\Msubspace{U}{k_1,..,k_p}{v_1,..,v_p}$, we obtain, using $\Psi$, all
%$\Amatrix{U'}{i'}{j'}$ spanning
%$\Msubspace{U'}{k'_1,..,k'_{p-1}}
%               {v_1,..,v_{q-1},v_{q+1},..,v_p}$.
\PROOFend %--------------------------------------------------------

Our most important result in this section is the following

%%%%%%%%%%%%%%%%%%%%%%%%%%%%%%%%%%%%%%%%%%%%%%%%%%%%%%%%%%%%%%%%%%%%%
%    THEOREM:   Algebraic sum is a direct sum
%%%%%%%%%%%%%%%%%%%%%%%%%%%%%%%%%%%%%%%%%%%%%%%%%%%%%%%%%%%%%%%%%%%%%

\begin{theorem}
  \label{theor_Mspace_as_direct_sum_of_M_subspaces}
  The algebraic sum of subspaces in
  (\ref{eq_Mspace_as_algebraic_sum_of_Msubspaces}) is a direct sum.
\end{theorem}

\PROOFstart %--------------------------------------------
We will prove the following statement equivalent to the above theorem.
\begin{quote}
  Let $\MsubspaceELEMENT{}{}{X}
                  \in \Msubspace{U}{}{}$,
      $\MsubspaceELEMENT{k_1,..,k_p}{v_1,..,v_p}{X}
                  \in \Msubspace{U}{k_1,..,k_p}{v_1,..,v_p}$
  for all possible choices of $({k_1,..,k_p})$ and $({v_1,..,v_p})$ be
  matrices whose sum is equal to a zero matrix. Then all these
  matrices (we will call them $X$-matrices) are zero matrices.
\end{quote}
The proof will be by induction with respect to $r$, the number of factors in
Kronecker product (\ref{eq_kronecker_product_of_unitaries}).

If we have only one factor, $U=\ithMATRIX{U}{1}$, then
$\Msubspace{U}{}{} = \Mspace{U}$ is a direct sum of one component,
namely $\Msubspace{U}{}{}$.

Assume now that the theorem is true for any $(r-1)$ factor Kronecker
product (\ref{eq_kronecker_product_of_unitaries}). In particular that what
we stated above as eguivalent to our theorem is true for such a
product.

Let our $X$-matrices add up to a zero matrix and let one of the
summands, $\MsubspaceELEMENT{k_1,..,k_p}{v_1,..,v_p}{X}$ be
nonzero. Let $-\MsubspaceELEMENT{k_1,..,k_p}{v_1,..,v_p}{X}$ be simply
denoted by $X$. Thus $X$ is the sum of the remaining $X$-matrices. We
also need to introduce $k'_1 < k'_2 < ... < k'_{r-p}$ such that
$\{k_1,..,k_p\} \cup \{k'_1,..,k'_{r-p}\} = \{1,..,r\}$.

The remaining $X$-matrices, forming $X$ in a sum, can be split into
two groups:
\begin{itemize}
  \item
        Into one group we put these belonging to subspaces
        $\Msubspace{U}{l_1,..,l_q}{w_1,..,w_q}$ such that
        $k'_1 \notin \{l_1,..,l_q\}$.
        ( An $X$-matrix belonging to $\Msubspace{U}{}{}$ is in this
        group, if $\{k_1,..,k_p\}$ is nonempty. The other situation we
        consider at the end of the proof.)
        Thanks to that choice we can
        apply Lemma \ref{lem_B_subrow_properties} \itemB\ to all of them, as
        well as to their sum, which will be denoted by $Y$:
        \begin{equation}
          \label{eq_Y_sum_subrow_properties}
          \sum_{c_{k'_1} = 1}^{n_{k'_1}}
            \ELEMENTof{Y}
                      {b}
                      {\vectorINDEX{:,..,:,c_{k'_1},:,..,:}}
          \ \ =\ \
          \ZEROvect\ \ \ \ \ \mbox{for any $b$-th row}\ \ .
        \end{equation}

  \item Into the second group we put all the $X$-matrices which belong to
      $\Msubspace{U}{l_1,..,l_q}{w_1,..,w_q}$ such that $k'_1 \in
      \{l_1,..,l_q\}$, that is $l_s = k'_1$ for some $s \in \{1,..,q\}$,
      and  then $w_s = v \in \{1,..,n_{k'_1}\}$. The group can further be
      split into subgroups corresponding to different values of $v$.
      When we add all the matrices in the whole group to form their sum $Z$,
      the $b$-th rows of $Z$,
      where the corresponding  $\vectorINDEX{b_1,..,b_r}$ satisfy $b_{k'_1}=v$,
      are sums of $b$-th rows
      of matrices belonging only to the subgroup associated with value $v$.
      Matrices from other subgroups have at these positions zero rows.

        We can apply  Lemma {\ref{lem_B_subrow_properties} \itemA} to
        each matrix in a subgroup, as well as to their sum $Z_{v}$:
        \begin{equation}
          \label{eq_Zv_sum_subrow_properties}
          \ELEMENTof{Z_v}
              {b}
              {\vectorINDEX{:,..,:,c_{k'_1},:,..,:}}
          \ \ =\ \
          \ABSOLUTEvalue{
            \ELEMENTof{\ithMATRIX{U}{k'_1}}{v}{c_{k'_1}}
          }^2
          \cdot
          \sum_{d_{k'_1} = 1}^{n_{k'_1}}
            \ELEMENTof{Z_v}
                      {b}
                      {\vectorINDEX{:,..,:,d_{k'_1},:,..,:}}
        \end{equation}
        Since, as we have just said, matrices forming $Z_v$ in a sum are
        responsible for formation of only those $b$-th rows for which
        $b_{k'_1} = v$, we can write for the sum $Z$ of all $Z_v$ :
        \begin{equation}
          \label{eq_Z_sum_subrow_properties}
          \ELEMENTof{Z}
                    {\vectorINDEX{b_1,..,b_{k'_1},..,b_r}}
                    {\vectorINDEX{:,..,:,c_{k'_1},:,..,:}}
          \ \ =\ \
          \ABSOLUTEvalue{
            \ELEMENTof{\ithMATRIX{U}{k'_1}}{b_{k'_1}}{c_{k'_1}}
          }^2
          \cdot
          \sum_{d_{k'_1} = 1}^{n_{k'_1}}
            \ELEMENTof{Z}
                      {\vectorINDEX{b_1,..,b_{k'_1},..,b_r}}
                      {\vectorINDEX{:,..,:,d_{k'_1},:,..,:}}
        \end{equation}

\end{itemize}

At this moment recall that, using our symbols, $X = Y + Z$. Also note that if
one had to, one would classify $X$ as belonging to the first group since $X
\in \Msubspace{U}{k_1,..,k_p}{v_1,..,v_p}$ and $k'_1 \notin \{k_1,..,k_p\}$.
So $X$ satisfies (\ref{eq_Y_sum_subrow_properties}) with $Y$ replaced by $X$
and analogous property holds for $Z$ through the linear combination $Z =
X-Y$. The last statement implies that the sums on the left hand side of
(\ref{eq_Z_sum_subrow_properties}) are zero rows, which in effect gives $Z =
\ZEROvect$. Further, let us return to component $Z_v$ of $Z$. As $Z$ is a
zero matrix and $Z_v$ have nonzero rows in disjoint regions, all $Z_v$ must
be zero matrices.

Recall that $Z_v$ are formed by $X$-matrices belonging to all possible
$\Msubspace{U}{l_1,..,l_q}{w_1,..,w_q}$ being subspaces of the larger space
$\MsubspaceFIXED{U}{k'_1}{v}$ introduced in Lemma
\ref{lem_isomorphic_Msubspaces} and spanned by
$\Amatrix{U}{i}{j}$, $\Smatrix{U}{i}{j}$ with
$i_{k'_1} = j_{k'_1} = v$.
Lemma \ref{lem_isomorphic_Msubspaces} tell us that this larger space
is isomorphic to $\Mspace{U'}$ where
$U' = \ithMATRIX{U}{1} \otimes \ithMATRIX{U}{k'_1 - 1} \otimes
      \ithMATRIX{U}{k'_1 + 1} \otimes \ithMATRIX{U}{r}$. We
denote this isomorphism by $\Psi$, that is
$\Psi\left( \MsubspaceFIXED{U}{k'_1}{v} \right) = \Mspace{U'}$.
In accordance with Lemma \ref{lem_isomorphic_Msubspaces} we choose
isomorphism $\Psi$ such that it maps any space
$\Msubspace{U}{l_1,..,l_q}{w_1,..,w_q}$ under consideration onto space
$\Msubspace{U'}{l'_1,..,l'_{q-1}}{w_1,..,w_{s-1},w_{s+1},..,w_q}$
where, as above, $l_s = k'_1,\ w_s = v,\ s \in \{1,..,q\}$,
and where
      $l'_1$, ..., $l'_{q-1}$
are positions of
      $l_1$ ..., $l_{s-1}$, $l_{s+1}$, ..., $l_q$
in the sequence
$(1,2,..., k'_1 - 1, k'_1 + 1, ..., r)$.
In particular,
  $\Msubspace{U}{l_s}{w_s} = \Msubspace{U}{k'_1}{v}$ is mapped onto
$\Msubspace{U'}{}{}$.

The images under $\Psi$ of $X$-matrices forming $Z_v$ are contained in the above
mentioned image spaces
$\Msubspace{U'}{l'_1,..,l'_{q-1}}{w_1,..,w_{s-1},w_{s+1},..,w_q}$, which in
turn are contained in $\Mspace{U'}$. The images add up to a zero matrix since the very
$X$-matrices forming $Z_v$ add up to a zero matrix.

On the other hand, and this is our induction assumption,
$\Mspace{U'}$ is a direct sum of
   $\Msubspace{U'}{l'_1,..,l'_{q-1}}{w_1,..,w_{s-1},w_{s+1},..,w_q}$,
for we have all possible parameterizing pairs of sequences coming from all
possible pairs $({l_1,..,l_q})$, $({w_1,..,w_q})$ such that $l_s = k'_1$,
$w_s = v$ for some $s \in \{1,..,q\}$.

( According to Def. \ref{def_Msubsets},
  sequences $({l_1,..,l_q})$ as well as $({l'_1,..,l'_{q-1}})$ are
  ordered.
  Let sequences $({l'_1,..,l'_{q-1}})$ and
  $({w_1,..,w_{s-1},w_{s+1},..,w_q})$ be given, with entries in
  appropriate ranges. Here $s$ is defined to be such that
  $l'_1 < \ldots < l'_{s-1} < k'_1$.
  Then for
  $({l_1,..,l_q})\ =\
  ({l'_1, \ldots, l'_{s-1}, k'_1, l'_s + 1, \ldots, l'_{q-1} + 1})$
  and
  $({w_1, \ldots, w_{s-1}, v, w_{s+1}, .., w_q})$ space
  $\Msubspace{U}{l_1,..,l_q}{w_1,..,w_q}$ is mapped onto
  $\Msubspace{U'}{l'_1,..,l'_{q-1}}{w_1,..,w_{s-1},w_{s+1},..,w_q}$. )

Thus all the images of $X$-matrices forming $Z_v$ are zero matrices and the
same must hold for the very $X$ matrices. As $v$ is taken to be an arbitrary
element from the set $\{1,..,n_{k'_1}\}$, all matrices forming $Z$, that is
belonging to the second group, are zero matrices. So $X$ is formed only by
matrices from the first group.

Analogous reasoning can be repeated with $k'_2$, ..., $k'_{r-p}$ instead of
$k'_1$. As a result we obtain that $X$ is formed only by $X$-matrices
belonging to $\Msubspace{U}{l_1,..,l_q}{w_1,..,w_q}$ such that $k'_1$, ...,
$k'_{r-p}$ do not belong to $\{l_1,..,l_q\}$ or, in other words, by $X$-matrices
belonging to $\Msubspace{U}{k_{g_1},..,k_{g_t}}{w_1,..,w_t}$, where
$\{k_{g_1},..,k_{g_t}\} \subset \{k_1,..,k_p\}$.

If one of these , call it $X'$, is nonzero and belongs to
$\Msubspace{U}{k_{g_1},..,k_{g_t}}{w_1,..,w_t}$ with $t < p$, that is
$\{k_{g_1},..,k_{g_t}\} \neq \{k_1,..,k_p\}$, it can be expressed with the
use of the remaining ones and $X$. We can repeat the whole reasoning with
$X'$ playing the role of $X$ and $k'_1$ being replaced by some $k' \in
\{1,..,r\} \setminus \{k_{g_1},..,k_{g_t}\}$ and simultaneously $k' \in
\{k_1,..,k_p\}$. We then classify summand $X$ as belonging to the second group of
matrices, since $X \in \Msubspace{U}{k_1,..,k_p}{v_1,..,v_p}$ and $k' \in
\{k_1,..,k_p\}$. We show as above that this second group of matrices
consists of zero matrices which contradicts our assumption that $X \neq
\ZEROvect$.

Note that at this stage $\Msubspace{U}{}{}$ could appear as
$\Msubspace{U}{k_{g_1},..,k_{g_t}}{w_1,..,w_t}$. The short version of our
reasoning for $X$ (the role of $X'$ at this stage) belonging to
$\Msubspace{U}{}{}$ is provided at the end of the proof. Also, precise
repetition of our reasoning is always possible because no components
in matrix sums are missing -- even if we write that a matrix is
formed only by some other matrices this means that the remaining ones
are there as zero matrices.

Having $X \in \Msubspace{U}{k_1,..,k_p}{v_1,..,v_p}$ on the left, on the
right we can only have nonzero components belonging to
$\Msubspace{U}{k_1,..,k_p}{w_1,..,w_p}$ with $(w_1,..,w_p) \neq
(v_1,..,v_p)$. Obviously $(w_1,..,w_p)$ are also distinct, since any
component $\Msubspace{U}{k_1,..,k_p}{w_1,..,w_p}$ in what we aim to prove to
be a direct sum, occurs only once. Collecting both sides together we get that
the sum of $X$-matrices ($-X$ is one of them) belonging to
$\Msubspace{U}{k_1,..,k_p}{w_1,..,w_p}$, for sequence $(k_1,..,k_p)$  fixed
by the choice of $X$ and for all possible sequences of values $(w_1,..,w_p)$,
is a zero matrix. Since regions where these last candidates for nonzero
$X$-matrices have nonzero rows do not overlap, also they must be zero
matrices. This is exactly what we wanted to prove, hence indeed
(\ref{eq_Mspace_as_algebraic_sum_of_Msubspaces}) is indeed the direct
sum.

Now, as promised twice in this proof, we consider the case of
$X \in \Msubspace{U}{}{}$, that is the case of $\{k_1,..k_p\}$ being
empty. Splitting the remaining summands into two groups, as
previously, for any $k' \in \{k'_1,..,k'_r\} = \{1,..,r\}$, we conclude
that any $X$-matrix belonging to
$\Msubspace{U}{l_1,..,l_q}{w_1,..,w_q}$ such that
$k' \in \{l_1,..,l_q\}$ is a zero matrix. Taking all $k'$ under
consideration we find that all $X$-matrices forming $X$ must be zero
matrices. So $X$ must be a zero matrix too.

\PROOFend %--------------------------------------------

Let us rewrite  (\ref{eq_Mspace_as_algebraic_sum_of_Msubspaces}) in a direct
sum fashion:
\begin{eqnarray}
  \label{eq_Mspace_as_direct_sum_of_Msubspaces}
  \lefteqn{\Mspace{U}\ \  =}  &  &  \\
%---------------------------------------------
  \nonumber
  &
  \Msubspace{U}{}{}
  \ \ \oplus \ \
  \bigoplus_{k \in \{1..r\}}
    \bigoplus_{v \in \{1..n_k\}}
      \Msubspace{U}{k}{v}
  \ \ \oplus\ \       &     \\
  %----------------------------------------------------------
  \nonumber
  &
  \ \ \oplus\ \
  \bigoplus_{k_1 < k_2 \in \{1..r\}}
    \bigoplus_{(v_1,v_2) \in \{1..n_{k_1}\} \times \{1..n_{k_2}\}}
      \Msubspace{U}{k_1,k_2}{v_1,v_2}
  \ \ \oplus \ \
  \ldots      &     \\
  %---------------------------------------------------
  \nonumber
  &
  \ldots \ \ \oplus \ \
  \bigoplus_{k_1 < .. < k_{r-1} \in \{1..r\}}
    \bigoplus_{(v_1,..,v_{r-1}) \in \{1..n_{k_1}\} \times .. \times \{1..n_{k_{r-1}}\}}
      \Msubspace{U}{k_1,..,k_{r-1}}{v_1,..,v_{r-1}}\ \ .
  &
\end{eqnarray}
Calculation of the defect of $U$ requires calculation of the dimension
of the $\Mspace{U}$, which luckily can be decomposed into the
dimensions of the direct sum components. We will make this
decomposition more informative in a moment.

%@@@@@@@@@@@@@@@@@@@@@@@@@@@@@@@@@@@@@@@@@@@@@@@@@@@@@@@@@@@@@@@@@@@@@@@@@@@@@@@@@@@@@@@@@@@@@@@
%
%                 E S T I M A T I N G    T H E    D E F E C T
%
%@@@@@@@@@@@@@@@@@@@@@@@@@@@@@@@@@@@@@@@@@@@@@@@@@@@@@@@@@@@@@@@@@@@@@@@@@@@@@@@@@@@@@@@@@@@@@@@

\section{Estimating the defect}
\label{sec_estimating_defect}

Calculation of all the dimensions of the direct sum components in
(\ref{eq_Mspace_as_direct_sum_of_Msubspaces}) can be reduced  to calculation
of dimensions of subspaces like $\Msubspace{U}{}{}$ considered for
subproducts of the Kronecker product $U$ in
(\ref{eq_kronecker_product_of_unitaries}).

%%%%%%%%%%%%%%%%%%%%%%%%%%%%%%%%%%%%%%%%%%%%%%%%%%%%%%%%%%%%
%     L E M M A :  On the dimension of M subspace
%%%%%%%%%%%%%%%%%%%%%%%%%%%%%%%%%%%%%%%%%%%%%%%%%%%%%%%%%%%%

\begin{lemma}
  \label{lem_dimension_of_Msubspace}
  There holds:
  \begin{equation}
    \label{eq_dimension_of_Msubspace}
    \dim\left(
      \Msubspace{U}{k_1,..,k_p}{v_1,..,v_p}
    \right)
    \ \ =\ \
    \dim\left(
      \Msubspace{U^\PRIMEstyle{p}}{}{}
    \right)
  \end{equation}
  where
  \begin{equation}
    \label{eq_reduced_kron_prod}
    U^\PRIMEstyle{p} \ \ =\ \
    \bigotimes_{k \in \{1..r\} \setminus \{k_1,..,k_p\}}
      \ithMATRIX{U}{k}
  \end{equation}
  is the original Kronecker product
  (\ref{eq_kronecker_product_of_unitaries})
  deprived of its $k_1$-th, ..., $k_p$-th factors.
\end{lemma}

\PROOFstart
In the proof we use Lemma \ref{lem_isomorphic_Msubspaces} repeatedly,
assigning appropriate values to $U$, $k$ and $v$. We have the
following relations:
\begin{itemize}
  \item Let $U = U$, the initial Kronecker product
      (\ref{eq_kronecker_product_of_unitaries}), and let $k = k_p$, $v =
      v_p$. Then by Lemma \ref{lem_isomorphic_Msubspaces}
      $\Msubspace{U}{k_1,..,k_p}{v_1,..,v_p}$ is isomorphic to
      $\Msubspace{U'}{k_1,..,k_{p-1}}{v_1,..,v_{p-1}}$, where $U' =
      \bigotimes_{l \in \{1..r\} \setminus \{k_p\}} \ithMATRIX{U}{l}$. In
      this case $k'_1$, ..., $k'_{p-1}$ mentioned in Lemma
      \ref{lem_isomorphic_Msubspaces} have values $k'_1=k_1$, ...,
      $k'_{p-1}=k_{p-1}$, as they are the positions of $k_1$, ...,
      $k_{p-1}$ in the sequence $(1,2,...,k_p -1, k_p +1,...,r)$.

  \item Let $U = U'$, $k = k_{p-1}$ and $v = v_{p-1}$. Then
      $\Msubspace{U'}{k_1,..,k_{p-1}}{v_1,..,v_{p-1}}$ is isomorphic to
      $\Msubspace{U''}{k_1,..,k_{p-2}}{v_1,..,v_{p-2}}$, where $U'' =
      \bigotimes_{l \in \{1..r\} \setminus \{k_{p-1}, k_p\}}
      \ithMATRIX{U}{l}$.

  \item
        ...

  \item
        Let $U = U^\PRIMEstyle{p-1} =
        \bigotimes_{l \in \{1..r\} \setminus \{k_2..k_p\}}
          \ithMATRIX{U}{l}$, and let $k = k_1$ and $v = v_1$. Then
        $\Msubspace{U^\PRIMEstyle{p-1}}{k_1}{v_1}$ is isomorphic to
        $\Msubspace{U^\PRIMEstyle{p}}{}{}$, where
        $U^\PRIMEstyle{p} =  \bigotimes_{l \in \{1..r\} \setminus \{k_1..k_p\}}
                    \ithMATRIX{U}{l}$.
\end{itemize}

Thus both subspaces mentioned in the lemma are isomorphic through the
composition of the above listed isomorphisms, so they are of equal
dimension.
\PROOFend %----------------------------------------------

\medskip
Let us simplify our notation by introducing, for $U$ being the product
(\ref{eq_kronecker_product_of_unitaries}):

\begin{equation}
  \label{eq_Msubspace_not_fixed_reduced_product_dim}
  \MsubspaceDIMENSION{U}{k_1,..,k_p}
  \ \ \stackrel{\mbox{def}}{=}\ \
  \dim\left( \Msubspace{U'}{}{} \right)
  \ \ \ \  \mbox{for}\ \ \ \
  U' =  \bigotimes_{k \in \{1..r\} \setminus \{k_1..k_p\}}
            \ithMATRIX{U}{k}
\end{equation}

\begin{equation}
  \label{eq_Msubspace_not_fixed_full_product_dim}
 \MsubspaceDIMENSION{U}{}
 \ \ \stackrel{\mbox{def}}{=}\ \
 \dim\left( \Msubspace{U}{}{} \right)
\end{equation}

Using the fact that $\Mspace{U}$ is the direct sum
(\ref{eq_Mspace_as_direct_sum_of_Msubspaces}) and that the dimension of any
component space $\Msubspace{U}{k_1,..,k_p}{v_1,..,v_p}$ can be calculated by
Lemma \ref{lem_dimension_of_Msubspace} to be
$\MsubspaceDIMENSION{U}{k_1,..,k_p}$, the dimension of $\Mspace{U}$ can be
expressed as in the Theorem \ref{theor_Mspace_dimension_formula} below.

At this point note also, that for a $(r-1)$-element set $\{k_1,...,k_{r-1}\}\
=\  \{1..r\} \setminus \{s\}$ there holds
$\MsubspaceDIMENSION{U}{k_1,..,k_{r-1}}\ =\
 \dim\left( \Msubspace{\ithMATRIX{U}{s}}{}{} \right)\ =\
 \dim\left( \Mspace{\ithMATRIX{U}{s}} \right)$. We do not assume any
 Kronecker product structure for $\ithMATRIX{U}{s}$, that is vector
 indices are identical with ordinary indices here. Using such vector
 indices we cannot fix their only component,
$\Msubset{\ithMATRIX{U}{s}}{1}{v}$ would be empty.

%therefore
%$\dim( \Mspace{\ithMATRIX{U}{s}})$ must be left as it is.

%%%%%%%%%%%%%%%%%%%%%%%%%%%%%%%%%%%%%%%%%%%%%%%%%%%%%%%%%%%%%%%%
%   T H E O R E M :  Formula for the dimencionality of
%                    the M space constructed for a Kronecker
%                    product of unitary matrices
%%%%%%%%%%%%%%%%%%%%%%%%%%%%%%%%%%%%%%%%%%%%%%%%%%%%%%%%%%%%%%%%
\begin{theorem}
  \label{theor_Mspace_dimension_formula}
  Let $U$ be a Kronecker product of unitary matrices
  (\ref{eq_kronecker_product_of_unitaries}). Then the dimension of the
  space $\Mspace{U}$  constructed for $U$ according to
  (\ref{eq_Mspace_definition}) reads:
  \begin{eqnarray}
    \label{eq_Mspace_dimension_formula}
    \lefteqn{
      \dim\left( \Mspace{U} \right)\ \ =
    }
    &    &    \\
%---------------------------------------------
    \nonumber
    &
    \MsubspaceDIMENSION{U}{}
    \ \ +\ \
    \sum_{k_1 \in \{1..r\}}
      n_{k_1} \cdot \MsubspaceDIMENSION{U}{k_1}
    \ \ +\ \
    \sum_{k_1 < k_2 \in \{1..r\}}
      n_{k_1} n_{k_2} \cdot \MsubspaceDIMENSION{U}{k_1,k_2}
    \ \ + \ldots
    &    \\
%-----------------------------------------------
    \nonumber
    &
    + \ \
    \sum_{k_1 < \ldots < k_{r-1} \in \{1..r\}}
      n_{k_1} \ldots n_{k_{r-1}} \cdot \MsubspaceDIMENSION{U}{k_1,..,k_{r-1}}
    &
  \end{eqnarray}
\end{theorem}
Accordingly (see Definition \ref{def_defect}), the defect of $U$ equals to
$(N-1)^2 \ -\  \dim\left( \Mspace{U} \right)$.

%On the applicability of formulas in case of 1 x 1 factors

Maybe it is a good moment to wonder about applicability of our formulas
when there are $1 \times 1$ Kronecker factors in product 
(\ref{eq_kronecker_product_of_unitaries}) forming $U$. 
This has been promised in the paragraph preceding Definition
\ref{def_Msubsets}.

First, recalling this definition, 
consider subset $\Msubset{U}{k_1,k_2,..,k_p}{v_1,v_2,..,v_p}$ of
set $\Mset{U}$. 
If  $k \in \{1,..,r\} \setminus \{k_1,..,k_p\}$ and $n_k = 1$, this
subset is empty. The reason is that you cannot find any
$\Amatrix{U}{i}{j}$, $\Smatrix{U}{i}{j}$ belonging to it, that is indexed
by $i,j$ satisfying $i_k \neq j_k$, 
               because $i_k,j_k \in \{1,..,n_k\} = \{1\}$ .
In this case let the corresponding 
              subspace $\Msubspace{U}{k_1,k_2,..,k_p}{v_1,v_2,..,v_p}$
of space $\Mspace{U}$ be defined as the zero space $\{ \ZEROvect \}$,
which is a space of dimension $0$. In particular, when there are $1
\times 1$ factors in  product
(\ref{eq_kronecker_product_of_unitaries}), 
$\Msubset{U}{}{}$ is empty and $\Msubspace{U}{}{}$ we define as the
zero space. Consequently, $\MsubspaceDIMENSION{U}{k_1,..,k_p}$ has value $0$,   
if  for some $k \in \{1,..,r\} \setminus \{k_1,..,k_p\}$ there holds
$n_k = 1$.
In other words $U'$, 
   used in the definition
   (\ref{eq_Msubspace_not_fixed_reduced_product_dim})
   of $\MsubspaceDIMENSION{U}{k_1,..,k_p}$,
must not contain any $1 \times 1$ Kronecker factors, if this 
$\MsubspaceDIMENSION{U}{k_1,..,k_p}$ is to be considered as a
potentially nonzero number. (It still can be zero if for example $U'$ is
a permutation matrix, in which case all $\Amatrix{U'}{i}{j}$,
$\Smatrix{U'}{i}{j}$ are zero matrices, so 
$\Msubspace{U'}{}{}$ is the zero space.)

Now, let us remove in sum (\ref{eq_Mspace_dimension_formula}) all
components $\MsubspaceDIMENSION{U}{k_1,..,k_p}$ that must be zero due
to the presence of $1 \times 1$ factors. What is left are only those 
$\MsubspaceDIMENSION{U}{k_1,..,k_p}$, for which there are no $1$'s 
left in sequence $(n_1,..,n_r)$ after removing the $k_1$-th, ..., $k_p$-th
entries. The sum (\ref{eq_Mspace_dimension_formula}) truncated in this
way looks as if it were written for $U$ deprived of its $1 \times 1$
Kronecker factors, in which case we accept the validity of formula 
(\ref{eq_Mspace_dimension_formula}). This new $U$ differs from the
original one by a unimodular factor, so the space $\Mspace{U}$ and the
defect $\DEFECT(U)$ remain the same, as one easily finds noting that
first of all $\Amatrix{U}{i}{j}$, $\Smatrix{U}{i}{j}$ do not change
with multiplication of $U$ by a unimodular number. The answer to our
question is thus: yes, we can use (\ref{eq_Mspace_dimension_formula})
in presence of $1 \times 1$ factors, provided that appropriate 
$\MsubspaceDIMENSION{U}{k_1,..,k_p}$ are assigned $0$'s. 

This has some practical importance. For example, if we wanted to
calculate defects for a large number of Kronecker products of at most
$r$ factors, we could store the sequences of sizes as rows of a matrix
with $r$ columns, filling rows with a certain number of $1$'s where
the number of factors would be less than $r$. These rows would then be
the input data for a procedure calculating the defect of an $r$
factor Kronecker product of unitaries.

We can go even furher. Let us also, in this context, return to the
direct sum formula (\ref{eq_Mspace_as_direct_sum_of_Msubspaces}). We
have adopted above a convention, according to which 
$\Msubspace{U}{k_1,k_2,..,k_p}{v_1,v_2,..,v_p}$ is the zero space if
there are $1$'s in sequence $(n_1,..,n_r)$ after removing from it the
$k_1$-th, ..., $k_p$-th entries. Let us throw out these zero spaces
from (\ref{eq_Mspace_as_direct_sum_of_Msubspaces}). We are left with 
$\Msubspace{U}{k_1,k_2,..,k_p}{v_1,v_2,..,v_p}$ where all $k$
associated with $n_k = 1$ sit in $\{k_1,..,k_p\}$. But every such
$\Msubspace{U}{k_1,k_2,..,k_p}{v_1,v_2,..,v_p}$ can be replaced by 
an equal space
$\Msubspace{U'}{k'_1,k'_2,..,k'_{p'}}{v'_1,v'_2,..,v'_{p'}}$, where:
\begin{itemize}
  \item
        $U'$ is $U$ deprived of its $1 \times 1$ Kronecker factors,
        with the order of the remaining factors preserved,
  \item
        $(k'_1,k'_2,..,k'_{p'})$ contains the positions, of subsequent
        entries of the subsequence obtained from $(k_1,..,k_p)$ by
        throwing out all $k$'s such that $n_k=1$, in the subsequence
        obtained from $(1,2,..,r)$ in the same way (This is
        cumbersome, being a consequence of our notation,
        and we met such formulation in Lemma \ref{lem_isomorphic_Msubspaces}.),
  \item
        %$(v'_1,v'_2,..,v'_{p'})$ is the result of throwing the
        %corresponding entries (of the only possible value $1$) out of 
        %$(v_1,v_2,..,v_p)$.
        $(v'_1,v'_2,..,v'_{p'})$ is the result of throwing out the 
        entries (of the only possible value $1$) out of 
        $(v_1,v_2,..,v_p)$ corresponding to entries $k$ of 
        $(k_1,..,k_p)$ for which $n_k = 1$.
\end{itemize}
(Both spaces are equal because they are spanned by the same matrices.
In general, $\Amatrix{U}{i}{j} = \Amatrix{U'}{i'}{j'}$,
where $i',j'$ correspond to the reduced vector indices
obtained from the vector indices corresponding to $i,j$, respectively, by removing those
$k$'th entries of these vector indices (entries necessarily equal to $1$), for which
$k$ satisfies $n_k = 1$.)

After this truncation and replacement 
(\ref{eq_Mspace_as_direct_sum_of_Msubspaces}) looks as if it were
written for $U'$, in which case we accept the validity of the direct
sum formula because $U'$ does not contain $1 \times 1$ Kronecker factors. 
Since $\Mspace{U} = \Mspace{U'}$, the only
difference between $U$ and $U'$ being a unimodular factor, the starting point
direct sum with all those zero spaces appears to lead to the correct
direct sum we would write if we first removed the $1 \times 1$
Kronecker factors. So, we can go either way, removing or leaving these
factors.

It is interesting to give the lower bound for the defect of $U$ being a
Kronecker product, because, as it was pointed in the Introduction, for most
unitaries the defect is zero. The lower bound we give here is associated with
upper bounds on the dimensions $\MsubspaceDIMENSION{U}{k_1,..,k_p}$ of
$\Msubspace{U'}{}{}$ for the Kronecker subproducts $U'$ as defined in
(\ref{eq_Msubspace_not_fixed_reduced_product_dim}), all of which are in fact
bounds on $\MsubspaceDIMENSION{U}{}$ with $U$ replaced by an appropriate
$U'$. All we need to know about these bounds is contained in the next
lemma. Note that although we assume that the sizes of
Kronecker factors are greater than $1$, any factor  of size $1$ -- a
unimodular number -- can be absorbed into one of the factors of size
greater than $1$.

%%%%%%%%%%%%%%%%%%%%%%%%%%%%%%%%%%%%%%%%%%%%%%%%%%%%%%%%%%%%%%%%%%%
%    L E M M A :  Bound on the dimensionalities of the M subspace
%                 of the M space, spanned by A^(i,j), S^(i,j) with
%                 i,j differeing at every subindex
%%%%%%%%%%%%%%%%%%%%%%%%%%%%%%%%%%%%%%%%%%%%%%%%%%%%%%%%%%%%%%%%%%%

\begin{lemma}
  \label{lem_Msubspace_not_fixed_dim_bounds}
  Let $U$ of size $N \times N$ be the Kronecker product
  (\ref{eq_kronecker_product_of_unitaries}) with $r \geq 1$ and
  factors $\ithMATRIX{U}{1}$, ..., $\ithMATRIX{U}{r}$ of size
  $n_1 \times n_1$, ..., $n_r \times n_r$, where all $n_k > 1$.

  Then the dimension of the space
  $\Msubspace{U}{}{}$ constructed for $U$ (Definition
  \ref{def_Msubspaces}) is bounded in one of the
  following ways:
  \begin{itemize}
    \item
          If $n_k > 2$ for any $k \in \{1..r\}$, then
          \begin{equation}
            \label{eq_bound_for_sizes_without_twos}
            \MsubspaceDIMENSION{U}{}
            \ \ \leq\ \
            (N-1)(n_1 - 1) \cdot \ldots \cdot (n_r -1)\ \ .
          \end{equation}

    \item
          If $n_{k_1} = \ldots = n_{k_p} = 2$ for some distinct values
          $k_s \in \{1..r\}$, and for $k \notin \{k_1..k_p\}$ there
          holds $n_k > 2$, then
          \begin{equation}
            \label{eq_bound_for_sizes_with_twos}
            \MsubspaceDIMENSION{U}{}
            \ \ \leq\ \
            (N-2^{p-1})(n_1 - 1) \cdot \ldots \cdot (n_r -1)\ \ .
          \end{equation}
  \end{itemize}
\end{lemma}

\PROOFstart
%OLD EXPLANATION
%Recall that $\Msubspace{U}{}{}$ is spanned by matrices
%$\Amatrix{U}{i}{j}$, $\Smatrix{U}{i}{j}$ such that the vectors indices
%corresponding to $i,j$ satisfy $i_1 \neq j_1$, ..., $i_r \neq
%j_r$. Thus, by Lemma \ref{lem_A_S_subrow_properties} \itemA,
%if $B \in \Msubspace{U}{}{}$ the subrows of any of its rows satisfy,
%for any $k \in \{1..r\}$ (where $c_k$ is at $k$-th position):
From Lemma \ref{lem_B_subrow_properties} \itemB, if $B \in \Msubspace{U}{}{}$
the subrows of any of its rows satisfy, for any $k \in \{1..r\}$ (where $c_k$
is at $k$-th position):
\begin{equation}
  \label{eq_B_subrow_equations}
  \sum_{c_k = 1}^{n_k}
    \ELEMENTof{B}{b}{\vectorINDEX{:,..,:,c_k,:,..,:}}
  \ \ =\ \
  \ZEROvect\ \ ,
\end{equation}
as in this case the  role of the set $\{k_1, ..., k_p\}$ from Lemma
\ref{lem_B_subrow_properties} is played by an empty set. We will show that
the $b$-th row of $B$ can be parametrized by at most
$(n_1 - 1) \cdot \ldots \cdot (n_r -1)$
independent parameters.
For example, let $\ELEMENTof{B}{b}{\vectorINDEX{c_1,...,c_r}}$ be known for
$c_1 \in \{1,..,(n_1 - 1)\}$, ..., $c_r \in \{1,..,(n_r - 1)\}$.
(\ref{eq_B_subrow_equations}) determines the remaining entries in the
following order:
\begin{itemize}
  \item
        $\ELEMENTof{B}{b}{\vectorINDEX{n_1,c_2,..,c_r}} =
          -\sum_{c_1 = 1}^{(n_1 - 1)}
             \ELEMENTof{B}{b}{\vectorINDEX{c_1,c_2,..,c_r}}$
        to make all $\ELEMENTof{B}{b}{\vectorINDEX{c_1,..,c_r}}$
        known for $c_1 \in \{1,..,n_1\}$,
                  $c_2 \in \{1,..,(n_2 - 1)\}$, ...,
                  $c_r \in \{1,..,(n_r - 1)\}$.
  \item
        $\ELEMENTof{B}{b}{\vectorINDEX{c_1,n_2,c_3,..,c_r}} =
          -\sum_{c_2 = 1}^{(n_2 - 1)}
             \ELEMENTof{B}{b}{\vectorINDEX{c_1,c_2,..,c_r}}$
        to make all $\ELEMENTof{B}{b}{\vectorINDEX{c_1,..,c_r}}$
        known for $c_1 \in \{1,..,n_1\}$,
                  $c_2 \in \{1,..,n_2\}$
                  $c_3 \in \{1,..,(n_3 - 1)\}$, ...,
                  $c_r \in \{1,..,(n_r - 1)\}$.
  \item
        ...
  \item
        $\ELEMENTof{B}{b}{\vectorINDEX{c_1,c_2,c_3,..,c_{r-1},n_r}} =
          -\sum_{c_r = 1}^{(n_r - 1)}
             \ELEMENTof{B}{b}{\vectorINDEX{c_1,c_2,..,c_{r-1},c_r}}$
        to have all the entries of $\ELEMENTof{B}{b}{:}$ known.
\end{itemize}
Next we show that the $b$-th row of $B$ filled using the above
algorithm, starting from those initially known
$(n_1 - 1) \cdot \ldots \cdot (n_r - 1)$ entries,
satisfies (\ref{eq_B_subrow_equations}) for any $k \in \{1..r\}$. That is,
that the following equality holds for any $k \in \{1..r\}$ and $c_s \in
\{1..n_s\}$ for any $s \in \{1..r\} \setminus \{k\}$:
\begin{equation}
  \label{eq_B_subrow_equations_other_version}
  \ELEMENTof{B}{b}{\vectorINDEX{c_1,..,c_{k-1},n_k,c_{k+1},..,c_r}}
  \ \ =\ \
  -\sum_{d_k = 1}^{n_k - 1}
     \ELEMENTof{B}{b}{\vectorINDEX{c_1,..,c_{k-1},d_k,c_{k+1},..,c_r}}
  \ \ .
\end{equation}
Let in the above expression  $c_{k_1} = n_{k_1}$, ..., $c_{k_p} =
n_{k_p}$ and $c_s \neq n_s$ for $s
\in \{1..r\} \setminus \{k_1,..,k_p\}$ with $k_1 < k_2 < \ldots <
k_p$. By this we mean that $n_k = c_k = c_{k_q} = n_{k_q}$ for some $q
\in \{1..p\}$. (Do not mistake $k_1$, ..., $k_p$ for those from
the second item of the lemma.)
Then the left hand side of
(\ref{eq_B_subrow_equations_other_version}), where some of the entries
are obtained using the algorithm from the dotted list above, reads:
\begin{eqnarray}
  \label{eq_B_entry_from_algorithm}
  \lefteqn{ \ELEMENTof{B}{b}{\vectorINDEX{c_1,..,c_r}} \ \ =}
  &  &  \\
%---------------------------------------------------------------
  \nonumber
  &
  -\sum_{d_{k_p} = 1}^{(n_{k_p} - 1)}
     \ELEMENTof{B}{b}{\vectorINDEX{c_1,..,d_{k_p},..,c_r}}
  \ \ =\ \
  -\sum_{d_{k_p} = 1}^{(n_{k_p} - 1)}
    \left(
      -\sum_{d_{k_{p-1}} = 1}^{(n_{k_{p-1}} - 1)}
         \ELEMENTof{B}{b}{\vectorINDEX{c_1,..,d_{k_{p-1}},..,d_{k_p},..,c_r}}
    \right)
  \ \ =\ \ \ldots
  &   \\
%---------------------------------------------------------------
  \nonumber
  &
  -\sum_{d_{k_p} = 1}^{(n_{k_p} - 1)}
     \left(
       -\sum_{d_{k_{p-1}} = 1}^{(n_{k_{p-1}} - 1)}
          \left(
            \ldots
            \left(
              -\sum_{d_{k_2} = 1}^{(n_{k_2} - 1)}
                 \left(
                   -\sum_{d_{k_1} = 1}^{(n_{k_1} - 1)}
                      \ELEMENTof{B}
                                {b}
                                {\vectorINDEX{c_1,..,d_{k_1},..,d_{k_p},..,c_r}}
                 \right)
            \right)
            \ldots
          \right)
     \right)
  \ \ =
  &    \\
%---------------------------------------------------------------
  \nonumber
  &
  -\sum_{d_{k_q} = 1}^{(n_{k_q} - 1)}
     \left(
       -\sum_{d_{k_p} = 1}^{(n_{k_p} - 1)}
          \left(
            ..%\ldots
            \left(
              -\sum_{d_{k_{q+1}} = 1}^{(n_{k_{q+1}} - 1)}
                 \left(
                   -\sum_{d_{k_{q-1}} = 1}^{(n_{k_{q-1}} - 1)}
                      \left(
                        ..%\ldots
                        \left(
                          -\sum_{d_{k_1} = 1}^{(n_{k_1} - 1)}
                             \ELEMENTof{B}
                                       {b}
                                       {\vectorINDEX{c_1,..,d_{k_1},..,d_{k_p},..,c_r}}
                        \right)
                        ..%\ldots
                      \right)
                 \right)
            \right)
            ..%\ldots
          \right)
     \right)
  &    \\
%-------------------------------------------------------------
  \nonumber
  &
  =\ \
  -\sum_{d_{k_q} = 1}^{n_{k_q} - 1}
     \ELEMENTof{B}
               {b}
               {\vectorINDEX{c_1,..,n_{k_1},..,n_{k_{q-1}},..,d_{k_q},..,n_{k_{q+1}},..,n_{k_p},..,c_r}}
  &   \\
%-------------------------------------------------------------
  \nonumber
  &
  =\ \
  -\sum_{d_{k_q} = 1}^{n_{k_q} - 1}
     \ELEMENTof{B}
               {b}
               {\vectorINDEX{c_1,..,c_{k_1},..,c_{k_{q-1}},..,d_{k_q},..,c_{k_{q+1}},..,c_{k_p},..,c_r}}
  &   \\
%-------------------------------------------------------------
  \nonumber
  &
  =\ \
  -\sum_{d_k = 1}^{n_k - 1}
     \ELEMENTof{B}
               {b}
               {\vectorINDEX{c_1,..,d_k,..,c_r}}\ \ ,
  &
\end{eqnarray}
and equals thus to the right hand side of
(\ref{eq_B_subrow_equations_other_version})

In this way each row of $B \in \Msubspace{U}{}{}$ could be potentially
parametrized by no more than $(n_1 - 1) \cdot \ldots \cdot (n_r - 1)$
parameters -- if there were no other restrictions caused by the structure of
$\ithMATRIX{U}{s}$ -- but we need to parameterize only $N-1$ rows in general.
This is caused by the property that all column sums in $B$ are zeros, just as
is the case with the spanning matrices $\Amatrix{U}{i}{j}$,
$\Smatrix{U}{i}{j}$. Hence, in general, we have no more than $(N-1) \cdot
(n_1 - 1) \cdot \ldots \cdot (n_r - 1)$ free parameters to determine
$B$. In other words: all $B \in \Msubspace{U}{}{}$ belong to a single
space of dimension $(N-1)(n_1 - 1) \cdot \ldots \cdot (n_r - 1)$.

The case when some $n_k$ are equal to $2$ is a separate one, but what really
counts is when there are more than one $n_k$ equal to $2$.

In this case let again $B \in \Msubspace{U}{}{}$ and let $n_{k_1} = 2$, ...,
$n_{k_p} = 2$. The $\vectorINDEX{b_1,..,b_r}$-th row of $B$ is a linear
combination of $\vectorINDEX{b_1,..,b_r}$-th rows of $\Amatrix{U}{i}{j}$,
$\Smatrix{U}{i}{j}$ such that, among other conditions, $i_{k_1} = 2/j_{k_1}$,
..., $i_{k_p} = 2/j_{k_p}$, where of course $i_{k_q}, j_{k_q} \in \{1,2\}$.

Consider the group of rows in $B$ indexed by such $\vectorINDEX{b_1,..,b_r}$
that for some chosen values $s_1,...,s_p \in \{1,2\}$ the vector index
satisfies
\begin{equation}
  \label{eq_index_scope_for_row_group_in_B}
  \left(
    b_{k_1} = s_1\   ,\  \ldots\ ,\ b_{k_p} = s_p
  \right)
   \ \ \ \ \mbox{or}\ \ \ \
   \left(
     b_{k_1} = 2/s_1\ ,\  \ldots\ ,\ b_{k_p} = 2/s_p
   \right)\ \ .
\end{equation}
Observe now that the rows of $B$ from this group are linear combinations of
respective nonzero rows of only those $\Amatrix{U}{i}{j}$,
$\Smatrix{U}{i}{j}$ (from the collection of those spanning
$\Msubspace{U}{}{}$) which are indexed by $i,j$ satisfying:
\begin{equation}
  \label{eq_AS_scope_for_spanning_row_group_in_B}
  \vectorINDEX{i_{k_1},..,i_{k_p}},\
  \vectorINDEX{j_{k_1},..,j_{k_p}}\ \
  \in\ \
  \{\ \vectorINDEX{s_1,..,s_p},\ \vectorINDEX{2/s_1,..,2/s_p}\ \}\ \ .
\end{equation}

Since subsets of $(i,j)$'s (of the set of all $(i,j)$'s indexing
$\Amatrix{U}{i}{j}$, $\Smatrix{U}{i}{j}$ spanning $\Msubspace{U}{}{}$)
fulfilling condition
(\ref{eq_AS_scope_for_spanning_row_group_in_B}) for fixed $s_1,...,s_p$ are
disjoint, also disjoint are the subsets of $\Msubset{U}{}{}$ of those corresponding
$\Amatrix{U}{i}{j}$, $\Smatrix{U}{i}{j}$ which are nonzero matrices.
   Within each such subset,
   associated with  (\ref{eq_AS_scope_for_spanning_row_group_in_B}),
   the nonzero rows of its members hit the index
   area defined in (\ref{eq_index_scope_for_row_group_in_B}).
Therefore the matrices $\Amatrix{U}{i}{j}$, $\Smatrix{U}{i}{j}$ belonging to
it not only have the property that their rows add up to a zero row, but also
that all the rows of such a matrix indexed by $\vectorINDEX{b_1,..,b_r}$
satisfying (\ref{eq_index_scope_for_row_group_in_B}) add up to a zero row.

Thus, through a linear combination, in the considered group of rows in $B$
one of the rows depends on others. We have $2^p/2 = 2^{p-1}$ possible choices
of $\vectorINDEX{s_1,..,s_p}$, hence $2^{p-1}$ groups (of the specified type)
of rows in $B$ with $2^{p-1}$ dependent rows altogether. Using the previous
method of parameterizing rows, we can parameterize no more than $(N-2^{p-1})$
rows, introducing no more than $(n_1 - 1) \cdot \ldots \cdot (n_r - 1)$
parameters into each row. In other words: $B$ being contained in the
described above $(N-2^{p-1})(n_1 - 1) \cdot \ldots \cdot (n_r - 1)$
dimensional space is a necessary condition for $B$ being a member of
$\Msubspace{U}{}{}$.

Our final remark will concern the case when we deal with a one factor
Kronecker product, $r=1$, $U = \ithMATRIX{U}{1}$, $N=n_1$.
The bound, either formula (\ref{eq_bound_for_sizes_without_twos})
or formula (\ref{eq_bound_for_sizes_with_twos}),
takes then the form $(N-1)^2$. This agrees with these facts:
\begin{itemize}
  \item
        $\dim(\Msubspace{U}{}{})\ =\ \dim(\Mspace{U})$, no Kronecker
        product structure assumed for $U$.
  \item
        $\Mspace{U}$ is contained in the space tangent to the set of
        all  $N \times N$ doubly stochastic matrices,
        that is the space of all  real $N \times N$ matrices in which
        entries in any row or column add up to $1$.
        This space is of dimension $(N-1)^2$.
\end{itemize}
\PROOFend %-------------------------------------------------------------

In what follows Theorem \ref{theor_Mspace_dimension_formula} there are
arguments behind defining $\Msubspace{U}{}{}$ as the zero space if
there are $1 \times 1$ Kronecker factors in
(\ref{eq_kronecker_product_of_unitaries}) not absorbed into larger
factors, and that $\MsubspaceDIMENSION{U}{}$ has to be assigned $0$
then. In this case we can also use the bound   
(\ref{eq_bound_for_sizes_without_twos}) or
(\ref{eq_bound_for_sizes_with_twos}) which produces the correct
dimension $0$.

Having the above result of Lemma 
\ref{lem_Msubspace_not_fixed_dim_bounds}
let us replace the values $\MsubspaceDIMENSION{U}{k_1,..,k_p}$ in
Theorem  \ref{theor_Mspace_dimension_formula}
by their upper bounds being the appropriate right hand sides of
(\ref{eq_bound_for_sizes_without_twos}) and
(\ref{eq_bound_for_sizes_with_twos}). We will get an upper bound on
$\dim\left( \Mspace{U} \right)$, equivalently a lower bound on the defect
$\DEFECT(U)\ =\ (N-1)^2 - \dim\left( \Mspace{U} \right)$. 

What is more, this can be done successfully also in the presence of
unabsorbed $1 \times 1$ Kronecker factors. Since the bounds
(\ref{eq_bound_for_sizes_without_twos}),
(\ref{eq_bound_for_sizes_with_twos})
yield $0$ for $\MsubspaceDIMENSION{U}{k_1,..,k_p}$ corresponding to
zero spaces in the direct sum
(\ref{eq_Mspace_as_direct_sum_of_Msubspaces}), whose presence is
associated with the $1 \times 1$ factors, bounding the dimensions of these spaces does
not affect the total upper bound on $\dim\left( \Mspace{U} \right)$ at
all, that is we obtain the value we would get if the $1 \times 1$
factors were first absorbed into larger factors. This can be better
understood if one carefully follows the reasoning following 
Theorem \ref{theor_Mspace_dimension_formula}. 
The practical consequence is this: 
if one had some other bounds on $\MsubspaceDIMENSION{U}{k_1,..,k_p}$, 
yielding zeros for spaces that are necessarily zero spaces because of
$1 \times 1$ factors, one could leave the $1 \times
1$ factors unabsorbed.

%But we donot stop here.

To better expose properties of our bound, it
is convenient to introduce another notion, the {\sl generalized defect}
$\GENdef(U)$ of a unitary $U$, defined as
\begin{equation}
  \label{eq_generalized_defect}    %%%%%%% DEF: Generalized defect
  \GENdef(U)\ \ =\ \
  \DEFECT(U) + (2N-1)\ \ =
  N^2 - \dim\left( \Mspace{U} \right)\ \ =
  \dim\left( \ORTHcomplement{\Mspace{U}} \right)
\end{equation}
When we recall expressions in
(\ref{eq_defect_calculation_with_orthogonal_complement}) and the
preceding formulas it will be clear that $\GENdef(U)$ is the dimension
of the space (\ref{eq_alowed_directions_parametrizing_space}), that is
the space:
\begin{equation}
  \label{eq_feasible_directions_space}
  \FEASIBLEparSPACE{U}\ \ =\ \
  \left\{
      R\ :\ \ \Ii R \HADprod U\ =\ E U\ \ \
      \mbox{for some antihermitian\ $E$}
  \right\}
\end{equation}

We use in our new definition the word {\sl generalized} because the definition of $\DEFECT(U)$ is
suited for the class of unitary matrices which have the dimension of
manifold
$\{ D_r \cdot U \cdot D_c\ :\ \ D_r,D_c\ \ \mbox{unitary diagonal}\}$
equal to $2N-1$. Because of some applications of the defect
mentioned in the Introduction it is more convenient subtract the dimension
of this manifold from $\GENdef(U)$ to define the defect associated
with other type of $U$.

The corresponding lower bound on $\GENdef(U)$ will be, thanks to
Theorem \ref{theor_Mspace_dimension_formula} and Lemma
\ref{lem_Msubspace_not_fixed_dim_bounds} and definitions
(\ref{eq_Msubspace_not_fixed_reduced_product_dim}),
(\ref{eq_Msubspace_not_fixed_full_product_dim}).
:
\begin{eqnarray}
  \nonumber
  (n_1 \cdot \ldots \cdot n_r)^2\ \ \ -\ \ \
  &  &      \\
%---------------------------------------------------
  \nonumber
  -\
  \left(
    \left(
      \prod_{l \in \{1..r\}}
        n_l
      \ -\ 2^{(\NUMBERofIN{\{1..r\}}{2}) - 1}
    \right)
    \cdot
    \prod_{l \in \{1..r\}}
      \left(n_l - 1 \right)
  \right.
  \ \ +
  &  &      \\
%---------------------------------------------------
  \nonumber
  \sum_{k_1 \in \{1..r\}}
    \left(
      n_{k_1}
      \left(
        \prod_{l \in \{1..r\} \setminus \{k_1\}}
          n_l
        \ -\ 2^{(\NUMBERofIN{\{1..r\} \setminus \{k_1\}}{2}) - 1}
      \right)
      \cdot
      \prod_{l \in \{1..r\} \setminus \{k_1\} }
        \left(n_l - 1 \right)
    \right)
  \ \ +
  &  &           \\
%----------------------------------------------------
  \nonumber
  \sum_{k_1 < k_2 \in \{1..r\}}
    \left(
      n_{k_1} n_{k_2}
      \left(
        \prod_{l \in \{1..r\} \setminus \{k_1,k_2\}}
          n_l
        \ -\ 2^{(\NUMBERofIN{\{1..r\} \setminus \{k_1,k_2\}}{2}) - 1}
      \right)
      \cdot
      \prod_{l \in \{1..r\} \setminus \{k_1,k_2\} }
        \left(n_l - 1 \right)
    \right)
  \ \ +\ \ \ldots\ \ +
  &  &           \\
%----------------------------------------------------
  \label{eq_generalized_defect_lower_bound}
  \left.
    \sum_{k_1 < ... < k_{r-1} \in \{1..r\}}
      \left(
        n_{k_1} ... n_{k_{r-1}}
        \left(
          \prod_{l \in \{1..r\} \setminus \{k_1..k_{r-1}\}}
            n_l
          \ -\ 2^{(\NUMBERofIN{\{1..r\} \setminus \{k_1..k_{r-1}\}}{2}) - 1}
        \right)
        \cdot
        \prod_{l \in \{1..r\} \setminus \{k_1..k_{r-1}\} }
          \left(n_l - 1 \right)
      \right)
  \right)\ \ ,
  &  &
\end{eqnarray}
where for a nonempty subset $\SET{A} = \{l_1,...,l_s\}$ of $\{1,...,r\}$ expression
$\NUMBERofIN{\SET{A}}{2}$ denotes
\begin{itemize}
  \item
        the number of $2$'s in sequence
        $\left( n_{l_1}, ..., n_{l_s} \right)$, if there are any,
  \item
        $1$ otherwise.
\end{itemize}

Although the above expression looks very complicated, below we show
that it is, almost,  a mere product of trivial functions of the sizes
$n_1$, ..., $n_r$ of factors of the considered Kronecker product.

Because in Lemma \ref{lem_Msubspace_not_fixed_dim_bounds} it was
assumed that all Kronecker factors are of size greater than $1$, formally
(\ref{eq_generalized_defect_lower_bound}) is a valid lower bound only
in this situation, therefore we confine ourselves in the formulation
of the theorem below. 
But, in the paragraph preceding the one 
containing the definition (\ref{eq_generalized_defect}) of the
generalized defect,
we argued that even in the presence of $1 \times 1$ Kronecker factors
the total bound (\ref{eq_generalized_defect_lower_bound}), obtained by
substituting 
(\ref{eq_bound_for_sizes_without_twos},
 \ref{eq_bound_for_sizes_with_twos}) into
(\ref{eq_Mspace_dimension_formula}),
will be equal to the value (\ref{eq_generalized_defect_lower_bound}) 
obtained for the Kronecker product with absorbed $1 \times 1$
factors 
(this can be also guessed purely by analyzing
(\ref{eq_generalized_defect_lower_bound}) as a function of sequence of
sizes $(n_1,..,n_r)$).
Let us leave it as a comment that
(\ref{eq_generalized_defect_lower_bound}) 
must, in the presence of $1 \times 1$ factors, lead to
(\ref{eq_GENdef_lower_bound_product_formula_without_twos}) or
(\ref{eq_GENdef_lower_bound_product_formula_with_twos} )
deprived of its factors $(2n_l - 1)$ corresponding to $n_l=1$, but
since then $2n_l -1\ =\ 1$, 
in Theorem \ref{theor_GENdef_lower_bound_product_formula} one does not
need to assume that all $n_k > 1$.

%%%%%%%%%%%%%%%%%%%%%%%%%%%%%%%%%%%%%%%%%%%%%%%%%%%%%%%%%%%%%%%
%
% THEOREM: The lower bound on the generalized defect D(U) of U
%          is the number:
%
%          a)   (2*n_1 -1) * ...... * (2*n_r - 1)
%                if all n_k are greater than 2
%
%          b)   (2*n_1 -1) * .... * (2*n_s - 1)*
%               2^((r-s)-1) * ( 2^(r-s) + 1 )
%               where n_1,...,n_r are all greater than 2,
%               if there are (r-s) 2's in the sequence
%               (n_1,...,n_r).  Here in fact it doesn't matter
%               that 2's are at the end of this sequence,
%               product (2*n_1 -1) * .... * (2*n_s - 1)
%               in fact is the product for all n_k > 2
%
%%%%%%%%%%%%%%%%%%%%%%%%%%%%%%%%%%%%%%%%%%%%%%%%%%%%%%%%%%%%%%%%

\begin{theorem}
  \label{theor_GENdef_lower_bound_product_formula}
  The lower bound
  (\ref{eq_generalized_defect_lower_bound}) on the generalized defect
  (definition (\ref{eq_generalized_defect}))
  of $U$ being the Kronecker product
  (\ref{eq_kronecker_product_of_unitaries}) with factors
  $\ithMATRIX{U}{1}$, ..., $\ithMATRIX{U}{r}$ of size $n_1 \times
  n_1$, ..., $n_r \times n_r$ respectively, where all $n_k > 1$, is equal to
  \begin{description}
    \item[\itemA]
          \begin{equation}
            \label{eq_GENdef_lower_bound_product_formula_without_twos}
            \prod_{l \in \{1..r\}}
              \left( 2 \cdot n_l\ -\ 1 \right),
          \end{equation}
          if all $n_k > 2$,

    \item[\itemB]
          \begin{equation}
            \label{eq_GENdef_lower_bound_product_formula_with_twos}
            \left(
              \prod_{l \in \{1..r\},\ n_l > 2}
                \left( 2 \cdot n_l\ -\ 1 \right)
            \right)
            \cdot
            2^{(\NUMBERofIN{\{1..r\}}{2}) - 1}
            \left(
              2^{(\NUMBERofIN{\{1..r\}}{2})}
              \ +\
              1
            \right)\ \ ,
          \end{equation}
          if there is at least one $2$ in the sequence
          $(n_1,...,n_r)$, where $\NUMBERofIN{\{1..r\}}{2}$ is the
          number of $2$'s in this sequence.
  \end{description}
\end{theorem}

\PROOFstart %---------------------------------------------------------
We need to show that there holds an equality between two values,
one being a relatively simple polynomial expression:
\begin{description}
  \item[\itemA]
        \begin{equation}
          \label{eq_simple_poly_difference_without_twos}
          \left(
            \prod_{l \in \{1..r\}}
               n_l
          \right)
          \cdot
          \left(
            \prod_{l \in \{1..r\}}
              \left(
                \left(n_l - 1 \right)
                \ +\
                1
              \right)
          \right)
          \ \ -\ \
          \prod_{l \in \{1..r\}}
            \left(
              n_l\ +
              \left( n_l - 1 \right)
            \right)
        \end{equation}
        if all $n_k > 2$,

  \item[\itemB]
        \begin{equation}
          \label{eq_simple_poly_difference_with_twos}
          \left(
            \prod_{l \in \{1..r\}}
              n_l
          \right)
          \cdot
          \left(
            \prod_{l \in \{1..r\}}
              \left(
                \left(n_l - 1 \right)
                \ +\
                1
              \right)
          \right)
          \ \ -\ \
          \left(
            \prod_{l \in \{1..r\},\ n_l > 2}
              \left(
                n_l\ +
                \left( n_l - 1 \right)
              \right)
          \right)
          \cdot
          2^{(\NUMBERofIN{\{1..r\}}{2}) - 1}
          \left(
            2^{(\NUMBERofIN{\{1..r\}}{2})}
            \ +\
            1
          \right),
        \end{equation}
        if $n_k = 2$ for some $k \in \{1,...,r\}$.
\end{description}
and the other being the long expression in the most outer bracket in
(\ref{eq_generalized_defect_lower_bound}). This last expression will
be simply referred to as (\ref{eq_generalized_defect_lower_bound})'.
We will compare the above mentioned two quantities component by
component.

We start from the $p$-th component of
(\ref{eq_generalized_defect_lower_bound})', by which we mean:
\begin{equation}
  \label{eq_pth_component_of_bound_prim}
  \sum_{k_1 < ... < k_p \in \{1..r\}}
    \left(
      n_{k_1} ... n_{k_p}
      \left(
        \prod_{l \in \{1..r\} \setminus \{k_1..k_p\}}
          n_l
        \ -\ 2^{(\NUMBERofIN{\{1..r\} \setminus \{k_1..k_p\}}{2}) - 1}
      \right)
      \cdot
      \prod_{l \in \{1..r\} \setminus \{k_1..k_p\} }
        \left(n_l - 1 \right)
    \right)\ \ .
\end{equation}

\begin{description}
  \item[Step 1]     %%%%%%%%%%%%%%%%%%%%%%%%%%%%%%%%%%%%%%  STEP 1
The, call it, left part being the result of taking only the left
product from the inner bracket in
(\ref{eq_pth_component_of_bound_prim}) is equal to:
\begin{equation}
  \label{eq_left_part_of_pth_component}
  \sum_{k_1 < ... < k_p \in \{1..r\}}
    \left(
      \prod_{l \in \{1..r\} }
        n_l
      \cdot
      \prod_{l \in \{1..r\} \setminus \{k_1..k_p\} }
        \left(n_l - 1 \right)
    \right)\ \ .
\end{equation}
The above sum can also be found in the expansion of the left product
of (\ref{eq_simple_poly_difference_without_twos})
or (\ref{eq_simple_poly_difference_with_twos}). Any component of the
above sum is the result of
choosing, in the right subproduct of the left product
of
(\ref{eq_simple_poly_difference_without_twos},
\ref{eq_simple_poly_difference_with_twos}),
$p$ brackets at
positions $k_1$, ..., $k_p$ from which $1$'s will be taken to be
multiplied with $(n_l - 1)$'s from the remaining $r-p$ brackets.
Thus every left part (\ref{eq_left_part_of_pth_component}) for
$p=1..(r-1)$ has its counterpart in the left product of
(\ref{eq_simple_poly_difference_without_twos},
\ref{eq_simple_poly_difference_with_twos}).

On the other hand, in
the expansion of this product there are two components, namely
$(n_1 \cdot \ldots \cdot n_r)(n_1-1) \cdot \ldots \cdot (n_r-1)$ and
$(n_1 \cdot \ldots \cdot n_r)$ so far unexplained. We easily see that
the first one is the left part of the $0$-th component of
(\ref{eq_generalized_defect_lower_bound})':
\begin{equation}
  \label{eq_0th_component_of_bound_prim}
  \left(
    \prod_{l \in \{1..r\}}
      n_l
    \ -\ 2^{(\NUMBERofIN{\{1..r\}}{2}) - 1}
  \right)
  \cdot
  \prod_{l \in \{1..r\}}
    \left(n_l - 1 \right)
\end{equation}
The second one still have to be found.

\item[Step 2]  %%%%%%%%%%%%%%%%%%%%%%%%%%%%%%%%%%%%%   STEP 2
The, call it, right part of (\ref{eq_pth_component_of_bound_prim})
looks like this:
\begin{equation}
  \label{eq_right_part_of_pth_component}
  -
  \sum_{k_1 < ... < k_p \in \{1..r\}}
    \left(
      n_{k_1} ... n_{k_p}
      \cdot
      2^{(\NUMBERofIN{\{1..r\} \setminus \{k_1..k_p\}}{2}) - 1}
      \cdot
      \prod_{l \in \{1..r\} \setminus \{k_1..k_p\} }
        \left(n_l - 1 \right)
    \right)\ \ .
\end{equation}
The expression under the sum, given for some $k_1,...,k_p$, can be further transformed into:
\begin{itemize}
  \item
        \begin{equation}
          \label{eq_right_part_summand_version_1}
          \left(
            \prod_{\{q:\ n_{k_q} \neq 2 \}}
              n_{k_q}
          \right)
          \cdot
          2^{(\NUMBERofIN{\{1..r\}}{2}) - 1}
          \cdot
          \prod_{l \in \{1..r\} \setminus \{k_1..k_p\},\ n_l \neq 2 }
            \left(n_l - 1 \right)
        \end{equation}
        if there are no $2$'s in $(n_1,...,n_r)$ or if there
        are $2$'s  in $(n_{l_1},...,n_{l_{r-p}})$ where
        $\{l_1,..,l_{r-p}\}\ =\ \{1,..,r\} \setminus \{k_1,..,k_p\}$.

  \item
        \begin{equation}
          \label{eq_right_part_summand_version_2}
          \left(
            \prod_{\{q:\ n_{k_q} \neq 2 \}}
              n_{k_q}
          \right)
          \cdot
          2^{(\NUMBERofIN{\{1..r\}}{2}) - 1}
          \cdot
          \left(
            \prod_{l \in \{1..r\} \setminus \{k_1..k_p\},\ n_l \neq 2 }
              \left(n_l - 1 \right)
          \right)
          \cdot 2
        \end{equation}
        if there are no $2$'s in $(n_{l_1},...,n_{l_{r-p}})$ and
        there are some in $(n_{k_1},..,n_{k_p})$, where
        $\{l_1,..,l_{r-p}\}$ is as above. Note that in this situation
        the additional $2$ at the end is caused by the fact that
        $\NUMBERofIN{ \{l_1,..,l_{r-p}\} }{2} = 1$ due to the
        definition of $\NUMBERofIN{}{2}$.

\end{itemize}
So, in the sum (\ref{eq_right_part_of_pth_component}) we add
expressions like
(\ref{eq_right_part_summand_version_1}) or
(\ref{eq_right_part_summand_version_2}) depending on how
$2$'s are scattered among the considered subsequences of
$(n_1,...,n_r)$. And we do this over all $p=1...(r-1)$ in the total
sum in (\ref{eq_generalized_defect_lower_bound})'.

Next we will add only those summands in the sum of all right
parts (\ref{eq_right_part_of_pth_component}) over $p=1...(r-1)$,
which are associated with sequences $(k_1,...,k_p)$ containing a given fixed
subsequence $(k_{q_1},...,k_{q_t}) = (s_1,...,s_t)$ such that $n_{s_1}$, ...,
$n_{s_t}$ are all greater than $2$,
and which (these $(k_1,..,k_p)$'s)
are such that for any
$k_j$ not belonging to this subsequence $n_{k_j} = 2$. The sizes $n_{s_1}$, ...,
$n_{s_t}$ play the role $n_{k_q}$ of
(\ref{eq_right_part_summand_version_1}) and
(\ref{eq_right_part_summand_version_2}). All sequences are assumed to
be increasing ones.
\begin{itemize}
  \item
        For $p=t$ we add nothing to build $(k_1,...,k_p)$ out of
        $(k_{q_1},...,k_{q_t})$, therefore we have only one such
        summand, here $q_1=1,...,q_t = p$. The number of summands is
        equal to $\binom{ \NUMBERofIN{\{1..r\}}{2} }{0}$.
        Formula (\ref{eq_right_part_summand_version_1}) is used for a summand.

  \item For $p=t+1$, the number of summands is equal to $0$ if there
      are no $2$'s in $(n_1,...,n_r)$, or $\binom{
      \NUMBERofIN{\{1..r\}}{2} }{1}$, that is the number of ways in
      which we can extend $(k_{q_1},...,k_{q_t})$ to $(k_1,...,k_p)$
      by such $k$ that $n_k = 2$. For a summand the formula
      (\ref{eq_right_part_summand_version_1}) is used if there are at
      least two $2$'s in $(n_1,..,n_r)$, that is if there is at least
      one $2$ outside $(n_{k_1},...,n_{k_p})$. Otherwise formula
      (\ref{eq_right_part_summand_version_2}) is used for a summand.

  \item
        For $p=t+2$ the number of summands is $0$ if there is only one
        $2$ in $(n_1,...,n_r)$, or
        $\binom{ \NUMBERofIN{\{1..r\}}{2} }{2}$ if there are at least two. If there
        are more than two $2$'s in $(n_1,...,n_r)$ formula
        (\ref{eq_right_part_summand_version_1}) is used for a summand,
        otherwise we use (\ref{eq_right_part_summand_version_2})

  \item
       ...

  \item
       If $p = t + \NUMBERofIN{\{1..r\}}{2}$ and $p \leq r-1$ the number of summands
       is either $0$ or
       $\binom{ \NUMBERofIN{\{1..r\}}{2} }{ \NUMBERofIN{\{1..r\}}{2} }$.
       In this case there are no more $2$'s in $(n_1,...,n_r)$ outside
       $(n_{k_1},...,n_{k_p})$, therefore we use formula
       (\ref{eq_right_part_summand_version_2}) for this single
       summand.
       (This all applies if there are $2$'s in $(n_1,...,n_r)$. If
       there are no $2$'s, we end with $p=t$.)

 \item
       If $p = r-1$ and $\NUMBERofIN{\{1..r\}}{2} = r - t$, then
       $p = t + \NUMBERofIN{\{1..r\}}{2} - 1$, the number of summands
       is either zero or it is equal to
       $\binom{ \NUMBERofIN{\{1..r\}}{2} }{ \NUMBERofIN{\{1..r\}}{2} - 1}$.
       Since there is one $2$ left in $(n_1,...,n_r)$ outside
       $(n_{k_1},...,n_{k_p})$ we use
       (\ref{eq_right_part_summand_version_1}) for all summands
       considered here.
       (This again applies if there are $2$'s in $(n_1,...,n_r)$.)

\end{itemize}
Having calculated the number of summands in each particular case, the
sum of all of them, with sign $-$ inherited from
(\ref{eq_right_part_of_pth_component}):
\begin{equation}
  \label{eq_right_parts_sum_subsum_for_fixed_nk_neq_2_subsequence}
  -
  \sum_{p=1}^{r-1}
    \sum_{\begin{array}{c}
            (k_1,...,k_p)\ \mbox{s.t.} \\
            (k_{q_1},..,k_{q_t}) = (s_1,...,s_t) \\
            \mbox{for some}\ q_1,..,q_t
          \end{array}
          }
      \left(
        n_{k_1} ... n_{k_p}
        \cdot
        2^{(\NUMBERofIN{\{1..r\} \setminus \{k_1..k_p\}}{2}) - 1}
        \cdot
        \prod_{l \in \{1..r\} \setminus \{k_1..k_p\} }
          \left(n_l - 1 \right)
      \right)
\end{equation}
takes the following forms.
\begin{itemize}
  \item %=========================  NO 2'S , A SUMMARY
        If there are no $2$'s in the sequence $(n_1,..,n_r)$
        \begin{equation}
          \label{eq_the_right_parts_sum_subsum_for_number_of_2s_being_zero}
          -
          \left(
            \prod_{i=1}^{t}
              n_{s_i}
          \right)
          \cdot
          2^{(\NUMBERofIN{\{1..r\}}{2}) - 1}
          \cdot
          \left(
            \prod_{l \in \{1..r\} \setminus \{s_1..s_t\},\ n_l \neq 2 }
              \left(n_l - 1 \right)
          \right)
          \cdot
          1\ \ ,
        \end{equation}
       (\ref{eq_right_parts_sum_subsum_for_fixed_nk_neq_2_subsequence})
       equals to
        \begin{equation}
          \label{eq_step_2_no_2s_final_result}
          -
          \left(
            \prod_{i=1}^{t}
              n_{s_i}
          \right)
          \cdot
          \prod_{l \in \{1..r\} \setminus \{s_1..s_t\} }
          \left(n_l - 1 \right)
          \ \ .
        \end{equation}
        In fact, in this case it is better to analyze directly
        (\ref{eq_right_part_of_pth_component}), where, by taking $p =
        1,...,r-1$, one finds that any summand has its counterpart in
        the right product of
        (\ref{eq_simple_poly_difference_without_twos}). Namely,
        (\ref{eq_step_2_no_2s_final_result}) for
        $(s_1,...,s_t)=(k_1,...,k_p)$ is the result of choosing $n_l$
        from brackets at positions $k_1$, ..., $k_p$ in the right
        product of (\ref{eq_simple_poly_difference_without_twos}), to
        be multiplied by $(n_l-1)$ from the remaining brackets. The
        additional summands in an expansion of the right product of
        (\ref{eq_simple_poly_difference_without_twos}) are the
        product of all $n_k$'s, $n_1 \cdot \ldots \cdot n_r$, which
        compensates for an identical product in an expansion of the
        left product of
        (\ref{eq_simple_poly_difference_without_twos}), we searched
        for this earlier at {\bf Step 1}, and finally the product of
        all $(n_k-1)$'s which in turn can be found in the $0$'th
        component (\ref{eq_0th_component_of_bound_prim}) of
        (\ref{eq_generalized_defect_lower_bound})', on the other side
        of the equality between
        (\ref{eq_simple_poly_difference_without_twos}) and
        (\ref{eq_generalized_defect_lower_bound})' being proved. Thus
        we are completely done with the case when there are no $2$'s
        in $(n_1,...,n_r)$.

  \item   %============================== WITH 2'S, A SUMMARY
        If the are some $2$'s in the sequence $(n_1,...,n_r)$:
        \begin{description}
          \item[Case $1$]   %.......................Case 1
                 If $t + \NUMBERofIN{\{1,...,r\}}{2}\ \leq r-1$
                 then
                 (\ref{eq_right_parts_sum_subsum_for_fixed_nk_neq_2_subsequence})
                 equals to
                \begin{eqnarray}
                  \nonumber
                  -
                  \left(
                    \prod_{i=1}^{t}
                      n_{s_i}
                  \right)
                  \cdot
                  2^{(\NUMBERofIN{\{1..r\}}{2}) - 1}
                  \cdot
                  \prod_{l \in \{1..r\} \setminus \{s_1..s_t\},\ n_l \neq 2 }
                    \left(n_l - 1 \right)
                  \cdot
                  &  &  \\
                %------------------------------------------------------------------
                  \label{eq_the_right_parts_sum_subsum_for_number_of_2s_plus_p_less_than_r}
                  \cdot
                  \left(
                    \binom{ \NUMBERofIN{\{1..r\}}{2} }{0}\ +\
                    \binom{ \NUMBERofIN{\{1..r\}}{2} }{1}\ +\
                    \binom{ \NUMBERofIN{\{1..r\}}{2} }{2}\ +\
                    \ldots\ +\
                    \binom{ \NUMBERofIN{\{1..r\}}{2} }{ \NUMBERofIN{\{1..r\}}{2} }\ +
                    1
                  \right)\ \ ,
                  &  &
                \end{eqnarray}
                In this case we can extend $(s_1,...,s_t)$ to
                $(k_1,...,k_p)$ in such a way that
                $(n_{k_1},...,n_{k_p})$ contains all $2$'s from
                $(n_1,...,n_r)$. The summand associated with such
                $(k_1,...,k_p)$ is expressed with the use of
                (\ref{eq_right_part_summand_version_2}), in which
                we have an additional factor $2$. That is why
                expression of type
                (\ref{eq_right_part_summand_version_1}) standing
                in the top part of
                (\ref{eq_the_right_parts_sum_subsum_for_number_of_2s_plus_p_less_than_r})
                has to be added twice, and hence the last $1$ in
                the rightmost bracket of
                (\ref{eq_the_right_parts_sum_subsum_for_number_of_2s_plus_p_less_than_r}).
                This expression finally takes the form:
                \begin{equation}
                  \label{eq_Case_1_result}
                  -
                  \left(
                    \prod_{i=1}^{t}
                      n_{s_i}
                  \right)
                  \cdot
                  2^{(\NUMBERofIN{\{1..r\}}{2}) - 1}
                  \cdot
                  \left(
                    \prod_{l \in \{1..r\} \setminus \{s_1..s_t\},\ n_l \neq 2 }
                      \left(n_l - 1 \right)
                  \right)
                  \cdot
                  \left(
                    2^{(\NUMBERofIN{\{1..r\}}{2})} + 1
                  \right)\ \ .
                \end{equation}

          \item[Case $2$]  %.............................case 2
                If $t + \NUMBERofIN{\{1,...,r\}}{2}\ = r$, that
                is $t = r - \NUMBERofIN{\{1,...,r\}}{2}$ we get
                for
                (\ref{eq_right_parts_sum_subsum_for_fixed_nk_neq_2_subsequence})
                \begin{eqnarray}
                  \nonumber
                  -
                  \left(
                    \prod_{i=1}^{t}
                      n_{s_i}
                  \right)
                  \cdot
                  2^{(\NUMBERofIN{\{1..r\}}{2}) - 1}
                  \cdot
                  \prod_{l \in \{1..r\} \setminus \{s_1..s_t\},\ n_l \neq 2 }
                    \left(n_l - 1 \right)
                  \cdot
                  &  &  \\
                %------------------------------------------------------------------
                  \label{eq_the_right_parts_sum_subsum_for_number_of_2s_plus_p_equal_r}
                  \cdot
                  \left(
                    \binom{ \NUMBERofIN{\{1..r\}}{2} }{0}\ +\
                    \binom{ \NUMBERofIN{\{1..r\}}{2} }{1}\ +\
                    \binom{ \NUMBERofIN{\{1..r\}}{2} }{2}\ +\
                    \ldots\ +\
                    \binom{ \NUMBERofIN{\{1..r\}}{2} }{
                      \NUMBERofIN{\{1..r\}}{2} - 1 }
                  \right)\ \ ,
                  &  &
                \end{eqnarray}
                 In other words, outside the subsequence
                 $(n_{s_1},...,n_{s_t})$ there are only $2$'s in
                 the sequence $(n_1,...,n_r)$. Every summand
                 associated with $(k_1,...,k_p)$ produced out of
                 such $(s_1,...,s_t)$ is expressed with the use
                 of (\ref{eq_right_part_summand_version_1}). Note
                 that in the case of only one $2$ in
                 $(n_1,...,n_r)$ the sum of binomials containes
                 only $\binom{1}{0}$. The sum
                 (\ref{eq_the_right_parts_sum_subsum_for_number_of_2s_plus_p_equal_r})
                 takes the final form:
                \begin{eqnarray}
                  \nonumber
                  -
                  \left(
                    \prod_{i=1}^{t}
                      n_{s_i}
                  \right)
                  \cdot
                  2^{(\NUMBERofIN{\{1..r\}}{2}) - 1}
                  \cdot
                  \left(
                    \prod_{l \in \{1..r\} \setminus \{s_1..s_t\},\ n_l \neq 2 }
                      \left(n_l - 1 \right)
                  \right)
                  \cdot
                  \left(
                    2^{(\NUMBERofIN{\{1..r\}}{2})} - 1
                  \right)\ \ =
                  &  &                 \\
                %---------------------------------------------------
                  \nonumber
                  -
                  \left(
                    \prod_{i=1}^{t}
                      n_{s_i}
                  \right)
                  \cdot
                  2^{(\NUMBERofIN{\{1..r\}}{2}) - 1}
                  \cdot
                  \left(
                    \prod_{l \in \{1..r\} \setminus \{s_1..s_t\},\ n_l \neq 2 }
                      \left(n_l - 1 \right)
                  \right)
                  \cdot
                  \left(
                    2^{(\NUMBERofIN{\{1..r\}}{2})} + 1
                  \right)\ \ +
                  &  &         \\
                %------------------------------------------------------
                  \label{eq_Case_2_result}
                  \left(
                    \prod_{i=1}^{t}
                      n_{s_i}
                  \right)
                  \cdot
                  2^{(\NUMBERofIN{\{1..r\}}{2})}
                  \cdot
                  \left(
                    \prod_{l \in \{1..r\} \setminus \{s_1..s_t\},\ n_l \neq 2 }
                      \left(n_l - 1 \right)
                  \right)
                  &  &
                \end{eqnarray}
        \end{description}

        The reader should realize, that the {\bf Case $2$} occures
        only once, when $(s_1,...,s_t)$ is of maximal length. Then we
        get an additional component, the second part of
        (\ref{eq_Case_2_result}). In general we add expressions like
        (\ref{eq_Case_1_result}), over all possible nonempty sequences
        $(s_1,...,s_t)$, that is over all $t$ ranging from $1$ to
        $r - (\NUMBERofIN{\{1..r\}}{2})$. So let us first concentrate
        on them.

        As in the case without $2$'s in $(n_1,...,n_r)$ we note that
        any expression like (\ref{eq_Case_1_result}) has its
        counterpart in the right product of
        (\ref{eq_simple_poly_difference_with_twos}). Again it is a
        matter of choosing brackets from which we take $n_l$'s, at
        the positions corresponding to the positions of $s_1,...,s_t$ in the
        increasing subsequence $(k:\ n_k > 2)$ of $(1,...,r)$. Note
        that in this process we also deal with the situation when all
        $n_k > 2$ are used in (\ref{eq_Case_1_result}), in fact in the
        'standard' part of (\ref{eq_Case_2_result})
        ($(s_1,...,s_t)$ of maximum length, {\bf
          Case 2} above), which corresponds to taking only $n_l$'s
        from brackets in (\ref{eq_simple_poly_difference_with_twos}).

        What remains in an expansion of
        (\ref{eq_simple_poly_difference_with_twos}) are the values:
        \begin{eqnarray}
          \label{eq_simple_poly_with_2s_remaining_expr_nl_minus_1s}
          -\ \
          \left(
            \prod_{l \in \{1..r\},\ n_l > 2}
              \left(
                n_l - 1
              \right)
          \right)
          \cdot
          2^{(\NUMBERofIN{\{1..r\}}{2}) - 1}
          \left(
            2^{(\NUMBERofIN{\{1..r\}}{2})}
            \ +\
            1
          \right)\ \ \ \ \ \ \mbox{and}
          &  &        \\
        %----------------------------------------------------
          \label{eq_simple_poly_with_2s_remaining_expr_all_nks}
          \prod_{l \in \{1..r\}}
            n_l,\ \
        \end{eqnarray}
        where the last product was mentioned as unexplained at the
        end of {\bf Step $1$}.

        On the other side, in
        (\ref{eq_generalized_defect_lower_bound})', we have still left
        the right part of the $0$-th component
        (\ref{eq_0th_component_of_bound_prim}) of
        (\ref{eq_generalized_defect_lower_bound})', which can be
        written as:
        \begin{equation}
          \label{eq_0th_component_of_bound_prim_another_version}
          -\
          \left(
            \prod_{l \in \{1..r\},\ n_l > 2}
              \left(n_l - 1 \right)
          \right)
          \cdot
          2^{(\NUMBERofIN{\{1..r\}}{2}) - 1}
          \ \ ,
        \end{equation}
        the additional component in (\ref{eq_Case_2_result}) occuring
        only for $(s_1,...,s_t)$ of maximum length, therefore written
        as
        \begin{equation}
          \label{eq_Case_2_result_additional_component}
          \left(
            \prod_{\{l:\ n_l > 2\}}
              n_l
          \right)
          \cdot
          2^{(\NUMBERofIN{\{1..r\}}{2})}\ \ ,
        \end{equation}
        and finally all summands in the left parts
        (\ref{eq_right_part_of_pth_component}) for $p=1,...,r-1$
        which are associated with sequences $(k_1,...,k_p)$ such that
        there are only $2$'s in $(n_{k_1},...,n_{k_p})$. We add these
        summands in {\bf Step 3} which completes the proof.

\end{itemize}

\item[Step 3]   %%%%%%%%%%%%%%%%%%%%%%%%%%%%%% STEP 3
The expression under the sum in
(\ref{eq_right_part_of_pth_component}), in case of
$(k_1,...,k_p)$ such that
$(n_{k_1},...,n_{k_p})\ =\ (2,...,2)$, takes the form
\begin{itemize}
  \item
        \begin{equation}
          \label{eq_right_part_summand_nki_eq_2_version_1}
          2^{(\NUMBERofIN{\{1..r\}}{2}) - 1}
          \cdot
          \prod_{\{l:\ n_l > 2\} }
            \left(n_l - 1 \right)
        \end{equation}
        if the number of $2$'s in $(n_1,...,n_r)$ is greater than $p$,

  \item
        \begin{equation}
          \label{eq_right_part_summand_nki_eq_2_version_2}
          2^{(\NUMBERofIN{\{1..r\}}{2}) - 1}
          \cdot
          \left(
            \prod_{\{l:\ n_l > 2\} }
              \left(n_l - 1 \right)
          \right)
          \cdot
          2
        \end{equation}
        if the number of $2$'s in $(n_1,...,n_r)$ is equal to $p$,
        $p\  =\  (\NUMBERofIN{\{1..r\}}{2})$. This is again caused by
        the definition of $\NUMBERofIN{}{2}$. In this case
        $2^{(\NUMBERofIN{\{1..r\} \setminus \{k_1..k_p\}}{2})}$
        is equal to $2$ while there are no $2$'s among $n_k$ indexed
        by $k \in \{1..r\} \setminus \{k_1..k_p\}$.

\end{itemize}

Now we calculate how many summands of this type can be found in
the sum of all the right parts
(\ref{eq_right_part_of_pth_component}). It is assumed that at
least one $2$ can be found in $(n_1,...,n_r)$. The other case
has been completely resolved in the paragraph following formula
(\ref{eq_step_2_no_2s_final_result}).
\begin{itemize}
  \item
        For $p=1$ the number of summands is equal to
        $\binom{\NUMBERofIN{\{1..r\}}{2}}{1}$, the number of
        ways we choose a one element subsequence, containing
        only $2$'s, from $(n_1,...,n_r)$. We use
        (\ref{eq_right_part_summand_nki_eq_2_version_2}) for a
        summand if there is only one $2$ in $(n_1,...,n_r)$,
        otherwise we use
        (\ref{eq_right_part_summand_nki_eq_2_version_1}).

  \item
        For $p=2$ the number of summands is either $0$ if there
        is only one $2$ in $(n_1,...,n_r)$, or
        $\binom{\NUMBERofIN{\{1..r\}}{2}}{2}$, we choose a two element
        subsequence this time.
        We use
        (\ref{eq_right_part_summand_nki_eq_2_version_2}) for each
        summand if there are only two $2$'s in $(n_1,...,n_r)$,
        or we use
        (\ref{eq_right_part_summand_nki_eq_2_version_1}) if
        there are at least three.

  \item
        ...

  \item
        For $p=\NUMBERofIN{\{1..r\}}{2}$ the number of summands
        is  $\binom{\NUMBERofIN{\{1..r\}}{2}}{\NUMBERofIN{\{1..r\}}{2}}$ and we of
        course use formula
        (\ref{eq_right_part_summand_nki_eq_2_version_2})
        for this single summand.
\end{itemize}

We add all the considered summands below, where the negative sign is
inherited from (\ref{eq_right_part_of_pth_component}). In the first sum
we still use the form of a summand originating from
(\ref{eq_right_part_of_pth_component}):
\begin{eqnarray}
  \nonumber
  -
  \sum_{p=1}^{r-1}
    \sum_{\begin{array}{c}
            (k_1,...,k_p)\ \mbox{s.t.} \\
            (n_{k_1},...,n_{k_p}) = (2,2,...,2) \\
          \end{array}
          }
      \left(
        n_{k_1} ... n_{k_p}
        \cdot
        2^{(\NUMBERofIN{\{1..r\} \setminus \{k_1..k_p\}}{2}) - 1}
        \cdot
        \prod_{l \in \{1..r\} \setminus \{k_1..k_p\} }
          \left(n_l - 1 \right)
      \right)\ \ =
  &  &            \\
%-------------------------------------------------------
  \label{eq_right_parts_sum_subsum_for_nk_only_2s_subsequence}
  -\
  2^{(\NUMBERofIN{\{1..r\}}{2}) - 1}
  \cdot
  \prod_{\{l:\ n_l > 2\} }
    \left(n_l - 1 \right)
  \cdot
  \left(
    \binom{\NUMBERofIN{\{1..r\}}{2}}{1}\ +\
    \binom{\NUMBERofIN{\{1..r\}}{2}}{2}\ +\
    \ldots\ +\
    \binom{\NUMBERofIN{\{1..r\}}{2}}{\NUMBERofIN{\{1..r\}}{2}}
    \ +\
    1
  \right)
\end{eqnarray}
where the last $1$ is the result of applying formula
(\ref{eq_right_part_summand_nki_eq_2_version_2}), with that additional
multiplication by $2$, to the single summand
for $p\ =\ \NUMBERofIN{\{1..r\}}{2}$. The final result is thus:
\begin{equation}
  \label{eq_step_3_final_result}
  -\
  \left(
    \prod_{\{l:\ n_l > 2\} }
      \left(n_l - 1 \right)
  \right)
  \cdot
  2^{(\NUMBERofIN{\{1..r\}}{2}) - 1}
  \cdot
  2^{(\NUMBERofIN{\{1..r\}}{2})}
\end{equation}

This is the end. What is left for us to say is that:
\begin{itemize}
  \item
        (\ref{eq_step_3_final_result}) together with
        (\ref{eq_0th_component_of_bound_prim_another_version}),
        both within
        (\ref{eq_generalized_defect_lower_bound})', compensate for
        (\ref{eq_simple_poly_with_2s_remaining_expr_nl_minus_1s})
        being part of
        (\ref{eq_simple_poly_difference_with_twos}).

  \item
        (\ref{eq_Case_2_result_additional_component}), a part
        of (\ref{eq_generalized_defect_lower_bound})', is equal to
        (\ref{eq_simple_poly_with_2s_remaining_expr_all_nks}) coming
        from
        (\ref{eq_simple_poly_difference_with_twos}).
\end{itemize}

\end{description}  %%%%% END OF STEP1,STEP2,STEP3
\PROOFend %-----------------------------------------------------

Our calculation of the bound
(\ref{eq_GENdef_lower_bound_product_formula_without_twos},
\ref{eq_GENdef_lower_bound_product_formula_with_twos})
on the generalized defect $\GENdef(U)$ of
$U$ given by (\ref{eq_kronecker_product_of_unitaries}) is based on the
formula for the dimension of $\Mspace{U}$ given by
Theorem \ref{theor_Mspace_dimension_formula}, which in turn is based
on the decomposition of $\Mspace{U}$ into direct sum components
provided by Theorem \ref{theor_Mspace_as_direct_sum_of_M_subspaces}.
As far as the case with no $2$'s in the sequence of sizes
$(n_1,...,n_r)$ is concerned, we can use a simpler reasoning leading
to the corresponding bound
(\ref{eq_GENdef_lower_bound_product_formula_without_twos}). To this
end let us consider the following lemma, in which we use
spaces defined in (\ref{eq_feasible_directions_space}).

%%%%%%%%%%%%%%%%%%%%%%%%%%%%%%%%%%%%%%%%%%%%%%%
%  LEMMA :  If R,S parametrize fasible directions for U,V
%           then (R x S) parametrizes a feasible direction
%           for U x V
%%%%%%%%%%%%%%%%%%%%%%%%%%%%%%%%%%%%%%%%%%%%%%%
\begin{lemma}
  \label{lem_R_S_feasible_pars_then_RxS_feasible_par}
  Let $U$ and $V$ be unitary matrices of size $N \times N$ and $M
  \times M$ respectively.
  Let $R \in \FEASIBLEparSPACE{U}$ and
      $S \in \FEASIBLEparSPACE{V}$, where $R,S$ are real matrices of
      the size identical with that of $U,V$, respectively.

  Then $R \otimes S\ \in\ \FEASIBLEparSPACE{U \otimes V}$.
\end{lemma}

\PROOFstart %--------------------------------------------------
From the assumptions there holds
\begin{equation}
  \label{eq_feasible_directions_for_factors}
  \Ii R \HADprod U\ =\ E U
  \ \ \ \ \mbox{and}\ \ \  \
  \Ii S \HADprod V\ =\ F V
\end{equation}
for some  antihermitian matrices $E$, $F$.
Then
\begin{equation}
  \label{eq_feasible_direction_for_product}
  (\Ii R \HADprod U) \otimes
  (\Ii S \HADprod V)
  \ \ =\ \
  -(R \otimes S) \HADprod (U \otimes V)
  \ \ =\ \
  EU \otimes FV
  \ \ =\ \
  (E \otimes F)(U \otimes V)\ \ ,
\end{equation}
from which we get
\begin{equation}
  \label{eq_product_feasible_direction}
  \Ii (R \otimes S) \HADprod (U \otimes V)
  \ \ =\ \
  (-\Ii E \otimes F) (U \otimes V)\ \ ,
\end{equation}
where $(-\Ii E \otimes F)$ is antihermitian because
$E \otimes F$ is hermitian.
So  $R \otimes S \in \FEASIBLEparSPACE{U \otimes V}$.
\PROOFend %--------------------------------------------------

Therefore, if
$\left( \ithMATRIX{R}{k}_i \right)_{i=1..\GENdef(\ithMATRIX{U}{k})}$
are bases for $\FEASIBLEparSPACE{\ithMATRIX{U}{k}}$, $k=1..r$, then
$\ithMATRIX{R}{1}_{i_1} \otimes \ldots \otimes \ithMATRIX{R}{r}_{i_r}$
form a set of
$\GENdef(\ithMATRIX{U}{1}) \cdot \ldots \cdot
\GENdef(\ithMATRIX{U}{r})$
independent vectors within
$\FEASIBLEparSPACE{
  \ithMATRIX{U}{1} \otimes
  \ldots \otimes
  \ithMATRIX{U}{r}}$. Hence:

%%%%%%%%%%%%%%%%%%%%%%%%%%%%%%%%%%%%%%%%%%%%%%
%  COROLLARY: Supermultiplicativity of the generalized
%             defect, with respect to kronecker products
%%%%%%%%%%%%%%%%%%%%%%%%%%%%%%%%%%%%%%%%%%%%%%
\begin{corollary}
  \label{cor_gen_def_super_multiplicativity}
    For a Kronecker product of unitaries
    (\ref{eq_kronecker_product_of_unitaries}) the generalized defect
    is supermultiplicative:
    \begin{equation}
      \label{eq_gen_def_super_multiplicativity}
      \GENdef\left(
        \ithMATRIX{U}{1} \otimes
        \ldots \otimes
        \ithMATRIX{U}{r}
      \right)
      \ \
      \geq
      \ \
      \GENdef(\ithMATRIX{U}{1}) \cdot
      \ldots \cdot
      \GENdef(\ithMATRIX{U}{r})
    \end{equation}
\end{corollary}

The above mentioned set of
$\ithMATRIX{R}{1}_{i_1} \otimes \ldots \otimes \ithMATRIX{R}{r}_{i_r}$
needs not to be a basis for
$\FEASIBLEparSPACE{
  \ithMATRIX{U}{1} \otimes
  \ldots \otimes
  \ithMATRIX{U}{r}}$.
For example, for the $2 \times 2$ unitary Fourier matrix:
\begin{equation}
  \label{eq_Fourier_example}
  F_2 \ \ =\ \
  \frac{1}{\sqrt{2}}
  \left[
    \begin{array}{cc}
      1 &  1 \\
      1 & -1
    \end{array}
  \right]
\end{equation}
we have that
$\GENdef(F_2 \otimes F_2)\ >\
\GENdef(F_2) \cdot \GENdef(F_2)$, as
$\GENdef(F_2 \otimes F_2)\ =\ 10$ and
$\GENdef(F_2)\ = 3$.

The space $\FEASIBLEparSPACE{\ithMATRIX{U}{k}}$ constructed for a
unitary $\ithMATRIX{U}{k}$ of size $n_k \times n_k$ contains at least
$2n_k -1 $ independent real matrices, denoted further by
$\ithMATRIX{R}{k}_i$, $i = 1..(2n_k - 1)$, which are all zero matrices except for the
$l$-th row or $m$-th column filled all with ones, where $l=1..n_k$ and
$m=2..n_k$. The corresponding directions
$\Ii \ithMATRIX{R}{k}_i \HADprod \ithMATRIX{U}{k}$ were mentioned at
(\ref{eq_phasing_directions}) . Thus
\begin{equation}
  \label{eq_Uk_gen_def_estimation}
  \GENdef\left( \ithMATRIX{U}{k} \right)
  \ \ \geq\ \
  2n_k - 1\ \ ,
\end{equation}
and consequently, by
Corollary \ref{cor_gen_def_super_multiplicativity},
\begin{equation}
  \label{eq_Uks_kron_prod__gen_def_estimation}
  \GENdef\left(
    \ithMATRIX{U}{1} \otimes
    \ldots \otimes
    \ithMATRIX{U}{r}
  \right)
  \ \ \geq\ \
  \left( 2 n_1 - 1 \right)
  \cdot \ldots \cdot
  \left( 2 n_r - 1 \right)\ \ .
\end{equation}
This is equal to
(\ref{eq_GENdef_lower_bound_product_formula_without_twos}), which is
also valid if there is only one $2$ in the sequence of sizes
$(n_1,...,n_r)$.

If there is more than one matrix of size $2 \times 2$ among
$\ithMATRIX{U}{k}$, in
(\ref{eq_Uks_kron_prod__gen_def_estimation}) we get a lower bound which is
smaller than the lower bound
(\ref{eq_GENdef_lower_bound_product_formula_with_twos}).
This is caused by the fact that, for $x$ denoting the number of $2$'s
in $(n_1,...,n_r)$, there holds an equality for natural $x>1$:
\begin{equation}
  \label{eq_U2s_bound_inequality}
  2^{x-1} \cdot \left( 2^x + 1 \right)
  \ \ >\ \
  (2 \cdot 2 - 1)^x\ \ ,
\end{equation}
which is equivalent to
\begin{equation}
  \label{eq_U2s_bound_inequality_2}
  4^x - 3^x\ \ >\ \ 3^x - 2^x\ \ ,
\end{equation}
which in turn is true because, after dividing both sides by
$(4-3)\ =\ (3-2)\ =\ 1$, on the left the summands are greater then on
the right. So, in this situation, unsplitting the $2 \times 2$
matrices one by one from
$\ithMATRIX{U}{1} \otimes  \ldots \otimes \ithMATRIX{U}{r}$ is not a
good strategy and it is better to write
(by Corollary \ref{cor_gen_def_super_multiplicativity} and by
(\ref{eq_Uk_gen_def_estimation})):
\begin{eqnarray}
  \nonumber
  \GENdef\left(
    \ithMATRIX{U}{1} \otimes
    \ldots \otimes
    \ithMATRIX{U}{r}
  \right)
  \ \ \geq\ \
  \left(
    \prod_{n_k > 2}
      \GENdef\left( \ithMATRIX{U}{k} \right)
  \right)
  \cdot
  \GENdef\left(
    \bigotimes_{n_k = 2}
      \ithMATRIX{U}{k}
  \right)\ \ \geq
  &    &     \\
%----------------------------------------------
  \label{eq_estimating_stategy_with_twos}
  \left(
    \prod_{n_k > 2}
      \left( 2 n_k - 1 \right)
  \right)
  \cdot
  \GENdef\left(
    \bigotimes_{n_k = 2}
      \ithMATRIX{U}{k}
  \right)\ \ .
  &     &
\end{eqnarray}

It would be interesting to conceive a set of $2^{x-1} (2^x + 1)$
independent matrices of some elegant structure within
$\FEASIBLEparSPACE{
  \bigotimes_{n_k = 2}
      \ithMATRIX{U}{k}
}$, where $x$ is the number of $2$'s in $(n_1,...,n_r)$,
so that we could write
(we can, but basing on Theorem \ref{theor_GENdef_lower_bound_product_formula}):
\begin{equation}
  \label{eq_U2s_gen_def_estimation}
  \GENdef\left(
    \bigotimes_{n_k = 2}
      \ithMATRIX{U}{k}
  \right)\ \ \geq\ \
  2^{x-1} \left( 2^x + 1 \right)\ \ .
\end{equation}
(We mean here a guess similar to the above choice of $(2n_k - 1)$
matrices $\ithMATRIX{R}{k}_i$ in 
$\FEASIBLEparSPACE{\ithMATRIX{U}{k}}$.)
This, together with (\ref{eq_estimating_stategy_with_twos}), would
provide the bound
(\ref{eq_GENdef_lower_bound_product_formula_with_twos} ) of
Theorem \ref{theor_GENdef_lower_bound_product_formula}, probably more directly
than with the use of Theorems
\ref{theor_Mspace_dimension_formula}
and
\ref{theor_Mspace_as_direct_sum_of_M_subspaces}.
Definitely, taking the $3^x$ matrices
$\bigotimes_{n_k = 2}
   \ithMATRIX{R}{k}_{i_k}$, where factors come from the set
\begin{equation}
  \label{eq_trivial_Rs_for_size_2}
  \left\{
    \ \
    \left[
      \begin{array}{cc}
        1 & 1 \\
        0 & 0
      \end{array}
    \right],
    \ \
    \left[
      \begin{array}{cc}
        0 & 0 \\
        1 & 1
      \end{array}
    \right],
    \ \
    \left[
      \begin{array}{cc}
        0 & 1 \\
        0 & 1
      \end{array}
    \right]
    \ \
  \right\}\ \ ,
\end{equation}
is, due to (\ref{eq_U2s_bound_inequality}),
not enough.

%@@@@@@@@@@@@@@@@@@@@@@@@@@@@@@@@@@@@@@@@@@@@@@@@@@@@@@@@@@@@@@@@@@@@@@@@@@@@@@@@@@@@@@@@@@@@@@@
%
%                 C O N C L U S I O N S
%
%@@@@@@@@@@@@@@@@@@@@@@@@@@@@@@@@@@@@@@@@@@@@@@@@@@@@@@@@@@@@@@@@@@@@@@@@@@@@@@@@@@@@@@@@@@@@@@@

\section{Conclusions}
  \label{sec_conclusions}

In this article we have achieved two things.

The first seems to be of great importance in numerical calculation of
the defect (Definition \ref{def_defect}) for large unitary $N \times
N$ matrices  $U$ with a Kronecker product structure. 
This involves calculation of the dimension of a certain
space $\Mspace{U}$  (defined in (\ref{eq_Mspace_definition})),
effectively in $\REALS^{N^2}$, 
spanned by the set $\Mset{U}$ (defined in (\ref{eq_Mset_definition})),
associated with $U$. It amounts to calculating the column rank of a
certain matrix, call it $M(U)$, built on the entries of $U$, and can
be performed for example by using rank/svd of MATLAB. When $U$ has a
Kronecker product structure, $\Mspace{U}$ can be split into a direct
sum with a large number of components 
(see Theorem \ref{theor_Mspace_as_direct_sum_of_M_subspaces}). So,
instead of applying the procedure to $M(U)$ we can apply it to a number
of its submatrices. This will make the whole calculation more reliable
and may prevent divergence in svd used in calculation of the rank.

The second thing we have achieved in this work is the calculation,
based on the above result, 
of the lower bound on the generalized defect $\GENdef(U)$ 
(defined in (\ref{eq_generalized_defect}), being the defect plus
$2N-1$) when $U$ is a Kronecker product. The generalized defect is
supermultiplicative with respect to Kronecker subproducts of $U$ 
(Corollary \ref{cor_gen_def_super_multiplicativity}) and this allows
us to trivially retrieve the lower bound of 
Theorem \ref{theor_GENdef_lower_bound_product_formula} when the number
of $2 \times 2$ Kronecker factors does not exceed $1$. In the other
case Theorem \ref{theor_GENdef_lower_bound_product_formula} gives a
better bound than that obtained using supermultiplicativity of
$\GENdef(U)$. We conjecture that, in either case, the lower bound is
attained by most matrices with a fixed Kronecker product structure
(fixed, up to an order, sequence of sizes of Kronecker factors).

All the formulas, the one expressing the direct sum forming
$\Mspace{U}$ (see (\ref{eq_Mspace_as_direct_sum_of_Msubspaces})), 
the resulting one expressing the dimension of $\Mspace{U}$ 
(see (\ref{eq_Mspace_dimension_formula})), 
and the one expressing the lower bound 
(see (\ref{eq_generalized_defect_lower_bound}),  
in a compact form 
(\ref{eq_GENdef_lower_bound_product_formula_without_twos}) or
(\ref{eq_GENdef_lower_bound_product_formula_with_twos})) are valid
also when there are $1 \times 1$ Kronecker factors among those forming $U$,
yielding the correct values we would get if the $1 \times 1$ factors
were absorbed into larger factors. This property allows, in
calculation of $\dim(\Mspace{U})$ or its bound, the use of a
procedure requiring a fixed number of factors, where the $1 \times 1$\
\ $[1]$'s can be taken as potentially missing factors to extend a shorter
Kronecker product.

The author believes that the first of the above mentioned achievements
will open the way to easier calculation of the defect for large
matrices $U$ with a Kronecker product structure. Such calculations may
be performed in order to asess whether a given $U$ gives rise to a
smooth family of inequivalent unitaries (i.e. not obtained one from
another by multiplying rows and columns by unimodular numbers) with
the moduli of entries fixed at the values sitting in 
$\left[ \ABSOLUTEvalue{U_{i,j}} \right]_{i,j=1..N}$, and how large its
dimension can be. This is a part of a question about unitary preimagies
$V$ to a doubly stochastic matrix $\left[ \ABSOLUTEvalue{U_{i,j}}^2
\right]_{i,j=1..N}$, that is about $V$ such that 
$\ABSOLUTEvalue{V_{i,j}}\ =\ \ABSOLUTEvalue{U_{i,j}}$ for all
$i,j$. This motivation is wider described in a prequel \cite{Defect}
to this paper. 
A special case, 
$\left[ \ABSOLUTEvalue{U_{i,j}}^2 \right]_{i,j=1..N}
\ =\ \left[ 1/N \right]_{i,j=1..N}$, the search for complex Hadamard
matrices, attracts even more interest, and it was dicussed in the
Introduction.

%@@@@@@@@@@@@@@@@@@@@@@@@@@@@@@@@@@@@@@@@@@@@@@@@@@@@@@@@@@@@@@@@@@@@@@@@@@@@@@@@@@@@@@@@@@@@@@@
%
%                 A C K N O W L D G E M E N T S 
%
%@@@@@@@@@@@@@@@@@@@@@@@@@@@@@@@@@@@@@@@@@@@@@@@@@@@@@@@@@@@@@@@@@@@@@@@@@@@@@@@@@@@@@@@@@@@@@@@

\section{Acknowledgements}
The author would like to thank Karol \.Zyczkowski and Marek Ku\'s from
Center for Theoretical Physics (Polish Academy of Sciences) in Warsaw,
for numerous discussions on the subject and reading this article.

Financial support by the grant number N N202 090239 of
Polish Ministry of Science and Higher Education is gratefully acknowledged.

%@@@@@@@@@@@@@@@@@@@@@@@@@@@@@@@@@@@@@@@@@@@@@@@@@@@@@@@@@@@@@@@@@@@@@@@@@@@@@@@@@@@@@@@@@@@@@@@@@@@@@
%
%     THE BIBLIOGRAPHY
%
%@@@@@@@@@@@@@@@@@@@@@@@@@@@@@@@@@@@@@@@@@@@@@@@@@@@@@@@@@@@@@@@@@@@@@@@@@@@@@@@@@@@@@@@@@@@@@@@@@@@@@

\end{document}